\documentclass[11pt]{amsart}
\usepackage{geometry}
\usepackage{multirow}
\usepackage{graphicx,xcolor,verbatim}
\usepackage{amssymb,amsfonts,amsmath}
\usepackage{graphics,epsfig}
\usepackage{graphicx}
\usepackage{latexsym}
\usepackage{algorithm}
\usepackage{algorithmic}

\newcommand{\be}{\begin{equation}}
\newcommand{\ee}{\end{equation}}

\newtheorem{remark}{Remark}
\newtheorem{theorem}{Theorem}

\usepackage{mathrsfs}
\usepackage{hyperref}
\hypersetup{
    colorlinks=true,
    linkcolor=blue,
    filecolor=magenta,      
    urlcolor=cyan,
    }
\usepackage{cite}
\DeclareMathOperator\supp{supp}
\DeclareMathOperator*{\argmin}{arg\,min}
\usepackage{bm}

\usepackage[normalem]{ulem}
\usepackage{amsaddr}
\usepackage{orcidlink}
\usepackage{fontawesome5}
\usepackage{subcaption}
\usepackage{subcaption}

\definecolor{myblue}{rgb}{0.21, 0.34, 0.74}
\definecolor{myred}{rgb}{0.79, 0.0, 0.09}
\definecolor{mygreen}{rgb}{0, 0.32, 0}

\usepackage{pifont}
\usepackage{enumitem}

\title[Data-Informed Mathematical Characterization of Absorption Properties]{Data-Informed Mathematical Characterization of Absorption Properties in Artificial and Natural Porous Materials}

\author{
Elishan C. Braun{\small $^{1}$},
Gabriella Bretti{\small $^{1,4}$}\orcidlink{0000-0001-5293-2115},
Melania Di Fazio{\small $^{2}$}\orcidlink{0000-0002-4472-5631},
Laura Medeghini{\small $^{2}$}\orcidlink{0000-0001-8568-8614},
Mario Pezzella{\small $^{3,4,*}$}\orcidlink{0000-0002-1869-945X}
}

\thanks{{\small $^{1}$} CNR National Research Council of Italy - Institute for Applied Mathematics ``M. Picone", Via dei Taurini, 19 - 00185 Rome, Italy}
\thanks{{\small $^{2}$} Sapienza, University of Rome - Department of Earth Sciences, Piazzale Aldo Moro 5 - 00185 Rome, Italy}
\thanks{{\small $^{3}$} CNR National Research Council of Italy - Institute for Applied Mathematics ``M. Picone", Via P. Castellino, 111 - 80131 Naples, Italy}
\thanks{{\small $^{4}$} Member of the Italian INdAM Research group GNCS}
\thanks{\textsuperscript{*} Corresponding author \faEnvelope[regular] \href{mailto:mario.pezzella@cnr.it}{\texttt{mario.pezzella@cnr.it}} 
}

\begin{document}
\begin{abstract}
    \small In this work, we characterize the water absorption properties of selected porous materials through a combined approach that integrates laboratory experiments and mathematical modeling. Specifically, experimental data from imbibition tests on marble, travertine, wackestone and mortar mock-ups are used to inform and validate the mathematical and simulation frameworks. First, a monotonicity-preserving fitting procedure is developed to preprocess the measurements, aiming to reduce noise and mitigate instrumental errors. The imbibition process is then simulated through a partial differential equation model, with parameters calibrated against rough and smoothed data. The proposed procedure appears particularly effective to characterize absorption properties of different materials and it represents a reliable tool for the study and preservation of cultural heritage. 
\end{abstract}
\maketitle

\section{Introduction}
Mathematical modeling has been increasingly employed to understand and predict the complex degradation processes affecting cultural heritage (we refer to \cite{MACH2019,MACH2021} for a comprehensive overview). Natural stones and artificial materials used in historical buildings and archaeological sites can be regarded as porous systems subject to different damaging processes due to the exposure to a variety of environmental and meteorological factors. 
Among all the causes, water infiltration is particularly relevant as it contributes to the formation of reactive flows (cf. \cite{volumemach,Chen_2006} and references therein) and triggers processes such as salt-crystallization \cite{BRACCIALE_Bretti_Alt}, penetration of air pollution \cite{Carbonation1, Sulphation2} and biocolonization \cite{charola}. For this reason, the investigation of the capillary absorption properties of building materials is of high interest \cite{Ass, bc, vg, Bretti24} and may offer valuable insights into degradation mechanisms.

The primary aim of this work is to provide a realistic characterization of the properties of commonly used building materials by integrating mathematical models with experimental measurements. For this purpose, we designed dedicated laboratory experiments from which we collected and analyzed data. Our approach relies on a data-informed Partial Differential Equation (PDE) model built upon the well-established Darcy's law \cite{bear,barenblatt}  that describes the water uptake into a porous medium, accounting for some material-specific physical quantities, such as the diffusion rate and the residual saturation level, which are calibrated against data.  Therefore, we address the inverse problem of determining the optimal model parameters by comparing the data with the simulation outcomes of a second order numerical method. For the parameters optimization, a multi-grid \cite{multigrid} particle swarm optimization algorithm is devised.  Furthermore, in order to reduce the influence of potential instrumental errors during the calibration phase, we propose a differential fitting technique that retains the information from the experiments while yielding a monotonic imbibition function, consistent with its physical interpretation. 

To the best of our knowledge, the scientific literature lacks similar contributions and detailed characterizations of the materials addressed in this work, especially for eco-friendly restoration mortars of new generation and wackestone. Indeed, these natural and artificial geomaterials were usually described from a minero-petrographic point of view, rarely taking a mathematical approach into consideration.

The paper is organized as follows.  A description of the specific building materials we analyzed and of the laboratory experiments we conducted is provided in Section \ref{sec:material}. The PDE model is presented in Section \ref{sec: Mathematical_Model}, where a quadratically convergent simulation method is provided, as well. Section \ref{sec: DATA_An} is devoted to data analysis, where the fitting technique for data preprocessing is introduced. The inverse problem pertaining to the optimal parameters selection is there formulated, as well. The calibration outcomes are presented and discussed in Section \ref{sec:Calibration}.  Conclusions and future perspectives of our research are  outlined in Section \ref{sec:concl}.
\section{Experimental setting and data}\label{sec:material}
\subsection{Description of materials used in the laboratory experiments}

\begin{figure}[ht]
    \centering
    \begin{subfigure}{0.32\textwidth}
        \centering
        \includegraphics[width=\textwidth]{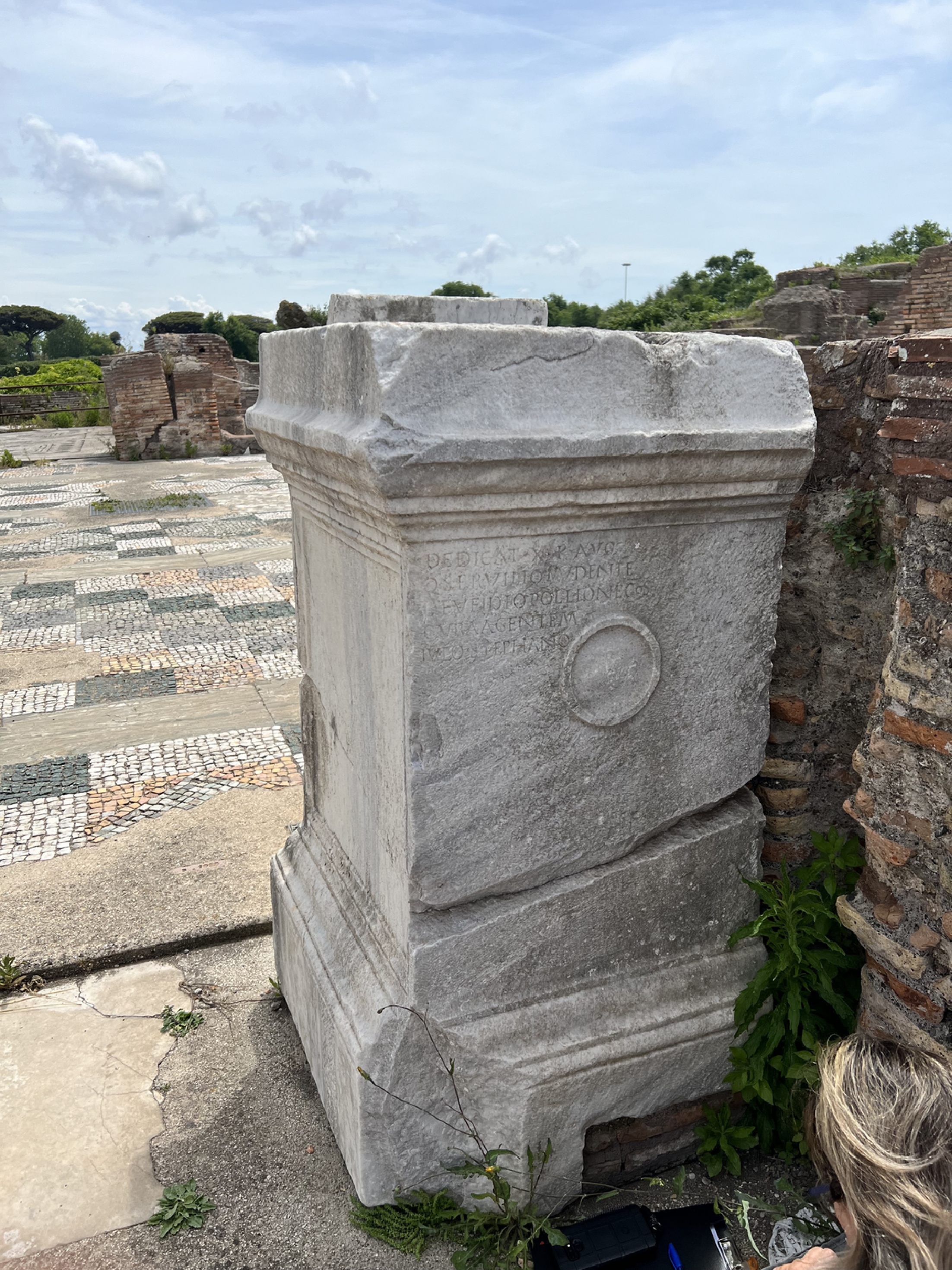}
        \caption{Marble altar in Ostia Antica, Lazio, Italy. \phantom{space here to align captions}}
        \label{fig:Altar_Marble}
    \end{subfigure}
    \hfill
    \begin{subfigure}{0.32\textwidth}
        \centering
        \includegraphics[width=\textwidth]{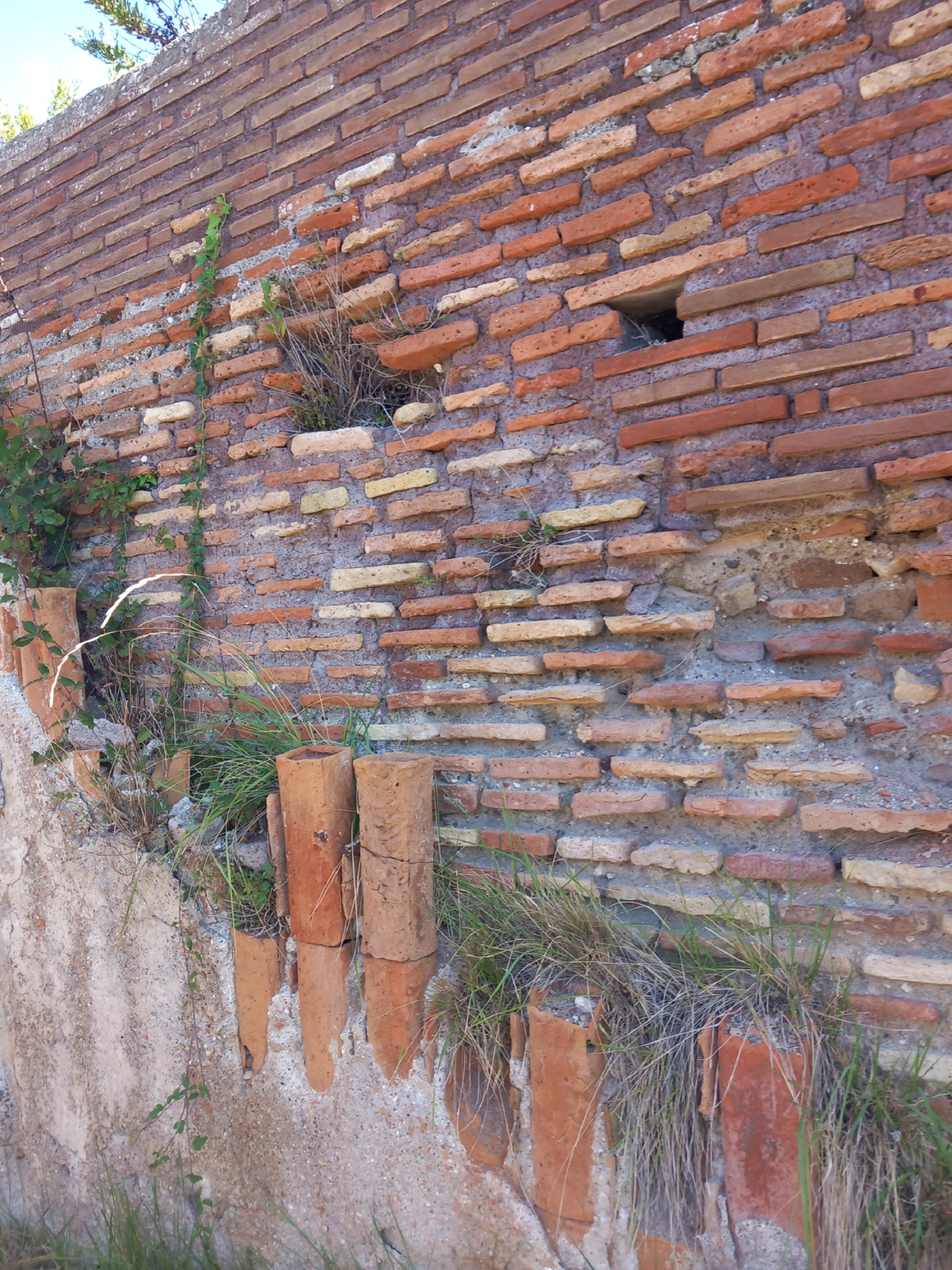}
        \caption{Mortar used in roman walls in Ostia Antica, Lazio, Italy.}
        \label{fig:Wall_Mortar}
    \end{subfigure}
    \hfill
    \begin{subfigure}{0.32\textwidth}
        \centering
        \includegraphics[width=\textwidth]{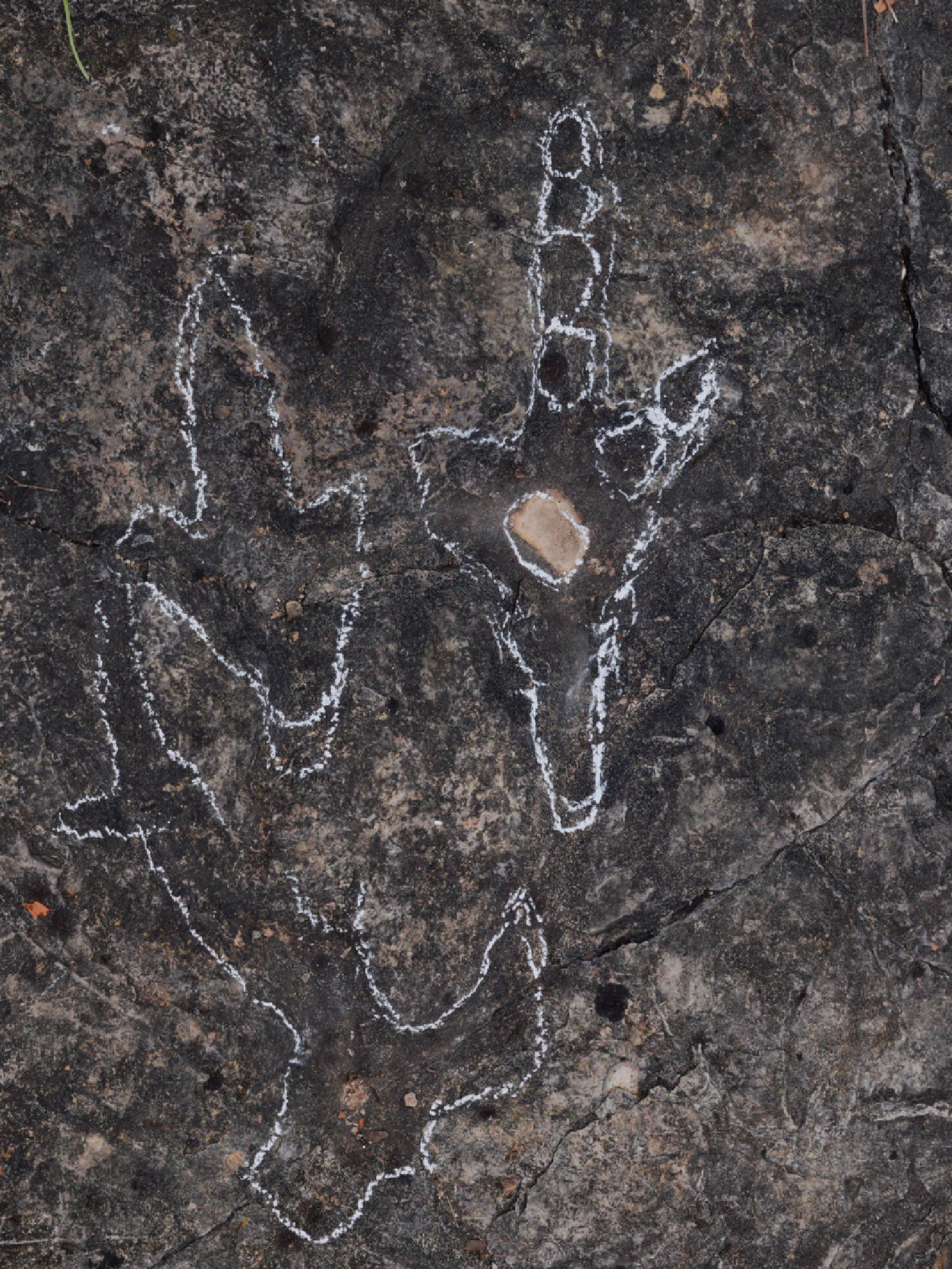}
        \caption{Wackestone  composing a natural monument in a former quarry in Sezze, Lazio, Italy.}
        \label{fig:Wackestone_Sezze}
    \end{subfigure}
    \caption{Examples of  natural stones and mortars in ancient buildings and archaeological sites. } 
    \label{fig:sites}
\end{figure}

Our analysis is performed through capillary water absorption tests conducted on samples of natural stones including white marble, travertine and wackestone, and artificial geomaterials such as hydraulic and aerial mortars. These materials are of interest due to their historical use in ancient constructions, with examples shown in Figure \ref{fig:sites}.
White marble is a metamorphic rock formed by the recrystallization of limestone or dolomite under high-temperature and high-pressure conditions. It is mainly composed of calcite crystals, with a granoblastic texture. Its porosity is usually very low and rarely exceeds 0.5\% of the mass (we refer to \cite{Gels,Malaga,Ruf,Zhang} for further details). Travertine is a chemically precipitated rock formed by the deposition of calcium carbonate from water rich in dissolved $\text{CO}_2$, typically from hydrothermal springs. Primarily consisting of calcite, it has a porous texture and exhibits color variations caused by mineral impurities. The travertine porosity is extremely variable and dependent on the deposit it comes from \cite{Mancini,Zalooli}. With the term wackestone is defined, according to the Dunham and Folk classification systems of limestones \cite{Adams_2017}, a fine-grained mud-supported carbonate sedimentary rock. It presents a fine grain size and micritic matrix, with a very poor porosity. Bioclasts and microfossils are present and distributed within the matrix.  

\begin{figure}[ht]
    \centering
    \begin{subfigure}{0.49\textwidth}
        \centering
        \includegraphics[scale=0.5]{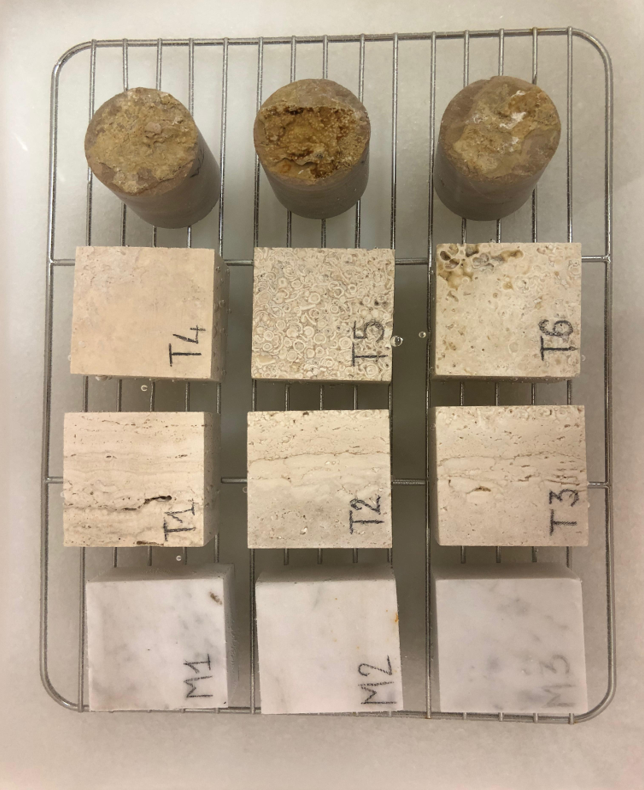}
        \caption{Natural stones specimens.}
        \label{fig:Materials_Natural}
    \end{subfigure}
    \begin{subfigure}{0.49\textwidth}
        \centering
        \includegraphics[scale=0.501]{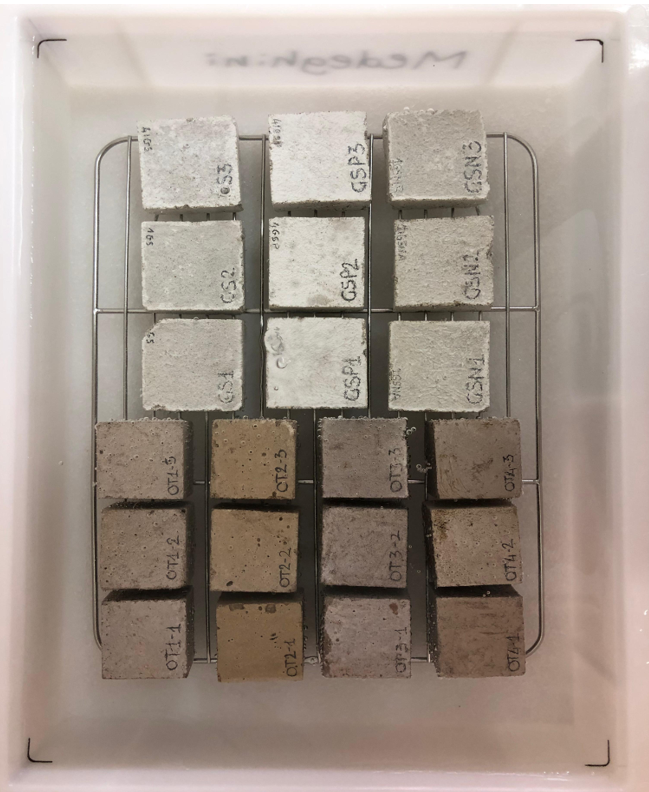}
        \caption{Artificial stones specimens.}
        \label{fig:Materials_Artificial}
    \end{subfigure}
    \caption{Materials samples employed for the laboratory capillary absorption tests.} 
    \label{fig:materials_new}
\end{figure}

For our experiments, we selected three samples of white Carrara marble and six samples of travertine from a cave in Tivoli, all in the form of cubes with an edge length of $5$ cm, as well as three wackestone specimens shaped as cylinders with a diameter of $4$ cm and a height of $5$ cm (see Figure \ref{fig:Materials_Natural}). Regarding the artificial geomaterials, we prepared two classes of specimens (cf. Figure \ref{fig:Materials_Artificial}). Twelve hydraulic mortars cubic samples (OT) of volume $64 \, \text{cm}^3$ were prepared based on an ancient Roman recipe, using natural hydraulic binder and pozzolan as aggregate \cite{Medeghini}. Furthermore, nine aerial mortar samples (GS, GSN and GSP) of size $4 \times 4 \times 2$ cm were produced by mixing lime and sand with the addition of organic additives extracted from microalgae and acrylic additives, in order to propose a green solution for the restoration \cite{Fratello2023}. We refer to Table \ref{table:allmaterials} for the ranges of porosity $\mathcal{N}$ of the monitored materials. It is important to note that for natural materials and aerial mortars, the porosity ranges are derived from scientific literature, as variations in composition can lead to different characteristics. In the case of hydraulic mortar samples instead, which consisted of four distinct compositions, porosity was experimentally determined by the producer through mercury intrusion porosimetry, yielding values of $33.5\%$ (OT1), $38.8\%$ (OT2), $40.2\%$ (OT3) and $44.3\%$ (OT4).

\begin{table}[htbp]
\centering
\begin{tabular}{|c|c|c|} \hline
\hline
 Material & Porosity range $\mathcal{N}$ & Reference  \\\hline\hline
Carrara Marble & $\phantom{1}0.50\% - \phantom{0}4.00\%$& \cite{Ruf}
\\
Tivoli Travertine & $\phantom{1}0.27\% - 16.67\%$& \cite{Mancini}
\\
Wackestone & $\phantom{1}0.05\% - \phantom{0}0.25\%$ &  -\\
GS-GSN-GSP & $13.00\% - 45.00\%$& \cite{Moropoulou,Santos,ZENG2012468}\\
OT1-OT2-OT3-OT4 & $33.50\% - 44.30\%$ & producer data\\
\hline\hline
\end{tabular}
\caption{Porosity ranges of monitored materials.}\label{table:allmaterials}
\end{table}

\subsection{Determination of capillary absorption coefficient and imbibition curve datapoints}
Water absorption by capillarity was evaluated according to UNI EN 1925:2000. Except for travertine, all samples exhibited uniform sedimentation planes, allowing an arbitrary choice of the absorption surface. For this specific material, however, tests were conducted on two sets of three samples: one with sedimentation planes parallel to the capillary rise direction and the other with planes perpendicular to it.

The tests were performed as follows. Samples of each material were preliminary dried in oven at $70 \pm 5 ^\circ C$, until a constant weight was attained, then stored in desiccator until they reached room temperature. Afterwards, they were placed in a plastic box with a lid and partially immersed in distilled water to a depth of $3 \pm 1 \ mm$.  Water absorption over time was assessed by periodically removing the samples from the water, drying their surface and then weighing them. 
Measurements were taken at the following time intervals: $1,\, 3,\, 5,\, 10,\, 15,\, 30,\, 60,\, 120,\, 180,\, 240$ minutes for artificial materials,  and up to 360 minutes for natural materials. After this initial phase, additional measurements were conducted every $24$ hours for nine days for the natural specimens, while for the artificial ones only two further measurements were taken at $24$ and $48$ hours.  Specifically, in compliance with the standards, we computed the amount of water absorbed per unit area by the specimens at time $t_i,$ i.e.
\begin{equation}\label{eq:Q_data}
    Q_{data}^i=\dfrac{m_i - m_d}{A}, \qquad \qquad i = 0, \dots, N_{data}, 
\end{equation}
expressed in $g/cm^{2},$ where 
\begin{itemize}
    \item $m_d$ is the mass of the dry specimen ($g$),
    \item $m_i$ is the mass of the specimen measured at time $t_i$ ($g$),
    \item $A$ is the area of the face immersed in water ($cm^2$).
\end{itemize}
We remark that the quantities in \eqref{eq:Q_data} were derived for each sample. However, for the sake of the analysis, we considered the average values computed across all the specimens of the same material.

\section{Mathematical and numerical frameworks}\label{sec: Mathematical_Model}
The starting point of this study is a one-dimensional mathematical model describing the imbibition process in a porous material characterized by a porosity $n_0$ and permeated by a liquid with density $\rho_l$ and viscosity $\mu_l$. Let $\theta_l(z,t)\in [0,n_0]$ represent the fraction of volume occupied by the fluid at height $z\in [0,\mathsf{H}]$ and time $t\in [0,\mathsf{T}].$ The mass balance equation then reads 
$\partial_t \theta_l(z,t)=\partial_z q(z,t),$ where $q$ is the volumetric flux expressed by Darcy's law (see, for instance \cite[Section 5.2]{bear} or \cite[Section 2.6]{barenblatt}) as follows
\begin{equation}\label{eq: Darcy}
q(z,t)=-\dfrac{k\!\left(\theta_l(z,t)/n_0\right)}{\mu_l} \, \left(\dfrac{\partial P_c}{\partial z} \! \!\left(\dfrac{\theta_l(z,t)}{n_0}\right) - \rho_l g\right), \qquad \quad (z,t)\in [0,\mathsf{H}] \times [0,\mathsf{T}].
\end{equation}

Here, $P_c(\cdot)$ denotes the capillary pressure, i.e. the pressure difference across the liquid-gas interface, $k(\cdot)$ is the intrinsic permeability of the porous matrix  to the vapor phase and $g$ represents the gravitational acceleration. 
 For the sake of simplicity, and for the entirety of this paper, we omit the subscript $l$ and denote $\theta_l(z,t)$ as $\theta(z,t).$ 
Following the arguments in \cite{clarelli2010,BRACCIALE_Bretti_Alt,Bretti24}, we introduce the fluid saturation $s(\theta(z,t))=\theta(z,t)/n_0\in [0,1]$ and the absorption function $B(s(\theta(z,t)))$ such that
\begin{equation}\label{eq: B_def}
\dfrac{\partial B}{\partial z}(s(\theta(z,t))) = -\dfrac{k(s(\theta(z,t)))}{\mu_l} \, \dfrac{\partial P_c}{\partial z}(s(\theta(z,t))), \qquad \qquad (z,t)\in [0,\mathsf{H}] \times [0,\mathsf{T}].
\end{equation}

The identification of the function $B(\cdot)$ that reproduces the capillary rise properties of a specific porous material is a challenging task. A straightforward and phenomenological approach is based on the observation that no fluid flows when the saturation is below the residual saturation $s_R\in[0,1)$, which ensures hydraulic continuity, or above the maximum saturation $s_S\in(s_R,1]$. Therefore, from a mathematical perspective, a reliable choice is to define $\frac{dB}{ds}(s)$ as a compactly supported function with $\supp(B^\prime(s))=[s_R, s_S]$ and a diffusion rate $D=\max_{[s_R, s_S]}B^\prime(s)$ that is calibrated against experimental data. We denote by $B^{\prime}(s)=\frac{dB}{ds}(s)$ and $B^{\prime \prime}(s)=\frac{d^2B}{ds^2}(s),$ for $s\in[0,1]$.
A third-degree polynomial form for $B(\cdot)$, which results in a concave parabola for its first derivative, is proposed in \cite{clarelli2010}. An alternative formulation, which accounts for the separate contributions of capillary pressure and permeability, is presented in \cite[Section 2.4.]{Bretti24}.

\begin{figure}[htp]
    \centering
    \includegraphics[width=0.70\linewidth]{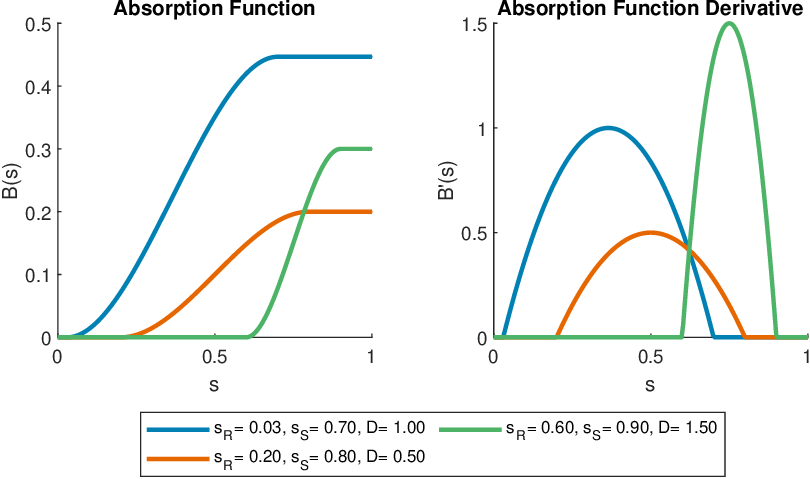}
    \caption{Absorption function and its derivative for different choices of parameters.}
    \label{fig:B_Plots}
\end{figure}

Here, in compliance with the empirical and phenomenological observations of \cite[Section 3.2]{clarelli2010}, we take 
\begin{equation}\label{eq: B'_shape}
    B^\prime(s)=\max\left\{0,-\frac{4 D (s_R - s) (s_S - s)}{(s_R - s_S)^2}\right\}, \qquad D=B^\prime\!\left(\frac{s_R+s_S}{2}\right),
\end{equation}
which yields
\begin{equation}\label{eq: B_shape}
    B(s)= 
    \begin{cases}
    0, & s \in [0,s_R],  \\  
    \dfrac{2D}{3} \, \dfrac{(s_R - s)^2 (-s_R + 3 s_S - 2 s)}{(s_R - s_S)^2}, \qquad \qquad & s \in [s_R,s_S], \\
     \dfrac{2D}{3}  \, (s_S - s_R), & s \in[s_S,1].
     \end{cases} 
\end{equation}
We refer to Figure \ref{fig:B_Plots} for some plots of the absorption function and its derivative for different choices of the parameters. 

The primary objective of this study is to determine the model parameters in \eqref{eq: B_shape} for the materials described in Section \ref{sec:material} by comparing numerical simulations and experiments. To this end, we define the overall quantity of fluid absorbed by the specimen at time $t$ as
\begin{equation}\label{eq: Q_defn}
    Q(t)=\rho_l \, \int_0^\mathsf{H} \! \theta(z,t) \ dz,
\end{equation}
which represents a physically measurable property of the system and therefore a crucial intermediary between the continuous PDE model \eqref{eq:richards} and the experimental data in \eqref{eq:Q_data}. Regarding the fraction of fluid $\theta,$ assuming small samples, we disregard the term $\rho_l g$ in \eqref{eq: Darcy} that accounts for gravitational effects. Hence, from \eqref{eq: Darcy} and \eqref{eq: B_def}, we retrieve the Richards equation modeling the penetration of water  
\begin{equation}\label{eq:richards}
  \dfrac{\partial \theta}{\partial t}(z,t) = \frac{\partial^2B}{\partial z^2}(s( \theta(z,t) ))=\frac{\partial}{\partial z}\left( \frac{1}{n_0} \dfrac{\partial B}{\partial s}( s( \theta(z,t) ) ) \, \dfrac{\partial\theta}{\partial z}(z,t)\right), 
\end{equation}
which is valid for $0 \leq z \leq \mathsf{H}$ and $0 \leq t \leq \mathsf{T}.$ To describe the onset of the imbibition process we set the initial condition 
\begin{equation}\label{eq:initial_condition}
    \quad \theta(h,0)=n_0 \chi_{\{0\}}(h)+\bar{\theta} \chi_{(0,\mathsf{H}]}(h), \qquad \qquad \qquad  h\in [0,\mathsf{H}],\tag{IC}
\end{equation}
where $\chi_{A}: A\to \{0,1\}$ denotes the characteristic function of the set $A\subset \mathbb{R}^+_0,$ while $\bar{\theta}$ describes the moisture content of the ambient air, that is assumed to be constant.
Then, denoted by $\texttt{SVD}$ the saturated vapor density, we consider
\begin{equation}\label{thetabar}
     \bar\theta = \texttt{SVD} \cdot \texttt{UR}  = \alpha_T \, (\alpha_0 + \alpha_1 T + \alpha_2 T^2 + \alpha_3 T^3) \cdot \texttt{UR} = 2.33 \cdot 10^{-5},
\end{equation}
where $T=27^\circ C$ is the average temperature during the experiment, $\texttt{UR}=80/100$ represents the percentage of humidity in the ambient air, $\alpha_T=10^{-6},$ $\alpha_0=5.02,$ $\alpha_1=0.32,$ $\alpha_2=8.18$ and $\alpha_3=3.12$ (we refer to \cite{HyperPhysics} for further details on the topic). 

We complement the PDE model \eqref{eq:richards} with the Dirichlet boundary conditions 
\begin{equation}\label{eq:DBC}
    \theta(0,t)=n_{0}, \qquad   \theta_l\left(\mathsf{H},t\right)=\bar{\theta}, \qquad \qquad  \ \qquad \qquad \ \ t\in [0,\mathsf{T}],\tag{DBC}
\end{equation}
or the Robin boundary conditions
\begin{equation}\label{eq:RBC}
    \theta(0,t)=n_{0}, \quad 
    \dfrac{\partial \theta}{\partial z}(\mathsf{H},t)= K_{w}(\bar{\theta}-\theta\left(\mathsf{H},t\right)), \qquad \quad \ t\in [0,\mathsf{T}],\tag{RBC}
\end{equation}
where $K_{w}\in \mathbb{R}^+$ is the water exchange rate. A numerical method for the simulation of problem \eqref{eq:richards} with the conditions \eqref{eq:DBC} and \eqref{eq:RBC} will be presented in Section \ref{sec: Numerical_Method}. 

\subsection{Numerical Method}\label{sec: Numerical_Method} 
In this section, we introduce a numerical method for approximating the solution of the equation \eqref{eq:richards}. Since our interest lies in solving the inverse problem of determining the model coefficients from data, it is crucial to opt for a reliable, efficient and straightforward integrator.

Consider the uniform spatial and temporal grids $z_j=j\Delta z$ and $t_k=k\Delta t,$ where $\Delta z$ and $\Delta t$ are the positive step-sizes. Let $N_z$ and $N_t$ be positive integers such that $\mathsf{H}=N_z \Delta z$ and $\mathsf{T}=N_t \Delta t.$ Denote by $\theta_j^k\approx \theta(z_j,t_k),$ for $j=0,\dots,N_z$ and $t=0,\dots,N_t.$ We propose the following Method of Lines (MOL) scheme 
\begin{equation}\label{eq: Numerical_Scheme}
    \begin{cases}
        \alpha_j^k=\left(B(\theta^k_{j+1}/n_0)-2B(\theta^k_j/n_0)+B(\theta^k_{j-1}/n_0)\right)/\Delta z^2, \\[0.17cm]
        \tilde{\theta}^k_j= \theta^k_j + \Delta t \, \alpha_j^k, \\[0.1cm]
        \beta_j^k=\left(B(\tilde{\theta}^k_{j+1}/n_0)-2B(\tilde{\theta}^k_j/n_0)+B(\tilde{\theta}^k_{j-1}/n_0)\right) /\Delta z^2, \\
        \theta_j^{k+1}=\theta^k_j + \frac{\Delta t}{2}(\alpha_j^k+\beta_j^k), \qquad \qquad \qquad \qquad \qquad \qquad \qquad \qquad \quad  \begin{array}{l}
           j=1,\dots,N_z-1,    \\
           k=0,\dots,N_t-1,   
        \end{array}
    \end{cases}
\end{equation}
derived by semi-discretizing the PDE with a second-order central finite difference for the spatial derivatives, followed by a second-order Heun’s approximation of the resulting ordinary differential equation (we refer to \cite{YUAN1999375,MIKHAIL198789} for futher details on MOL procedures). The initial condition \eqref{eq:initial_condition} and the boundary conditions \eqref{eq:DBC}-\eqref{eq:RBC} are discretized as follows
\begin{align}
    &\theta_j^0 = n_0 \chi_{\{0\}}(z_j) + \bar{\theta} \chi_{(0,\mathsf{H}]}(z_j), \qquad \qquad \qquad \qquad \qquad \qquad j = 0, \dots, N_z, \tag{dDIC}\label{eq:D_initial_condition} \\
    &\theta_0^k = n_0, \quad \quad \theta_{N_z}^k = \bar{\theta}, \qquad \qquad \qquad \qquad \qquad \qquad \qquad \quad \, \,   k = 0, \dots, N_t, \tag{dDBC}\label{eq:D_DBC} \\
    &\theta_0^k = n_0, \quad \quad \theta_{N_z}^k = \dfrac{4\theta_{N_z-1}^k - \theta_{N_z-2}^k + 2K_a \bar{\theta} \Delta z}{3 + 2K_a \Delta z},  \qquad \qquad k = 0, \dots, N_t. \tag{dDRBC}\label{eq:D_RBC}
\end{align}
Additionally, we approximate the observable in \eqref{eq: Q_defn} via a second order composite trapezoidal quadrature rule, i.e.
\begin{equation}\label{eq: Q_num}
    Q^k=\rho_l \, \frac{\Delta z}{2}\left(\theta_{0}^{k}+2\sum_{j=1}^{N_{x}}\theta_{j}^{k}+\theta_{N_{x}+1}^{k}\right)\approx Q(t_k), \qquad \qquad \qquad \quad k=0,\dots,N_t.
\end{equation}

The explicit MOL scheme \eqref{eq: Numerical_Scheme} is second-order accurate both in time and space, offering an optimal balance between computational cost and accuracy. This is particularly advantageous in contexts where model calibration based on experimental data is required, as such tasks require multiple iterations of the scheme. In such cases, maintaining an explicit form ensures that the resulting algorithm remains straightforward to implement and not computationally demanding. A standard Von Neumann analysis (see, for instance \cite{Crank_1996,Shishkina2007,Wesseling} and references therein) establishes the stability of the scheme \eqref{eq: Numerical_Scheme} when the CFL-like condition
\begin{equation}\label{eq: Stability_Condition}
    \Delta t\leq \frac{n_{0} \Delta z^{2}}{\max_{[0,1]}|B^\prime(s)|},
\end{equation}
is satisfied. 

\begin{figure}[htp]
    \centering
    \includegraphics[width=0.75\linewidth]{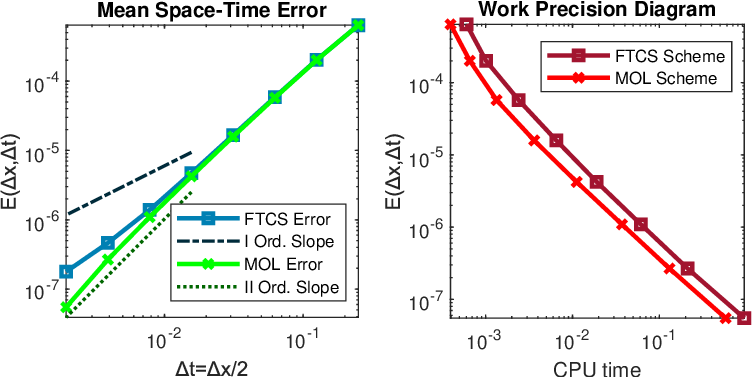}
    \caption{Experimental convergence and performances of the FTCS and MOL methods.}
    \label{fig:MOLvsFTCS}
\end{figure}

In order to experimentally highlight the advantages of the MOL discretization \eqref{eq: Numerical_Scheme}, we compare it with the FTCS (Forward in Time, Central in Space) scheme adopted in \cite[Section 3]{Bretti24}. Specifically, we approximate the solution of problem \eqref{eq:richards}-\eqref{eq:initial_condition}-\eqref{eq:DBC} with 
\begin{equation*}
    n_0=0.285, \quad \bar{\theta}=6.254 \cdot 10^{-2}, \quad  \mathsf{H}=8, \quad \mathsf{T}=60, \quad s_R=0.219, \quad s_S=1, \quad D=9.807\cdot 10^{-4}
\end{equation*}
and compare the mean space-time error and the experimental orders of convergence, defined as follows
\begin{equation}\label{eq:Error}
        E(\Delta z, \Delta t)=\dfrac{1}{N_zN_t}\sum_{j=0}^{N_z}\sum_{k=0}^{N_t}\left|\theta_j^k-\Theta_j^k\right|, \qquad \rho=\log_2\left(\dfrac{E(\Delta z, \Delta t)}{E(\frac{\Delta z}{2}, \frac{\Delta t}{2})}\right).
\end{equation}
Here, the reference solution $\Theta_j^k$ is computed by the MOL method with $\Delta z=2^{-9}$ and $\Delta t=2^{-12}$. Figure \ref{fig:MOLvsFTCS} and Table \ref{tab:Exp_Convergence} report the simulations outcomes for different values of the stepsizes and confirm the quadratic convergence of the MOL integrator \eqref{eq: Numerical_Scheme}. Furthermore, the work precision diagram in Figure \ref{fig:MOLvsFTCS} shows the mean errors with respect to the mean computational efforts over ten runs and confirms the superior performances of the MOL discretization \eqref{eq: Numerical_Scheme} compared to the FTCS approach in \cite[Section 3]{Bretti24}.

\begin{table}[htbp]
\begin{center}
  \begin{tabular}{|c|cc|cc|} \hline  \hline
   \multicolumn{1}{|c|}{ Stepsizes} & \multicolumn{2}{c|}{ Mean Space-time Errors} & \multicolumn{2}{c|}{ Exp. order}  \\
    $\Delta t=\Delta x/2$ &  $E_{FTCS}$ &  $E_{MOL}$ &  $\rho_{FTCS}$ &  $\rho_{MOL}$ \\ \hline 
    \hline \rule{-1.5pt}{9pt}
    $\!2^{-2\phantom{1}}$& $6.40\cdot 10^{-4}$  & $6.40\cdot 10^{-4}$ & $-$    & $-$    \\
    $2^{-3\phantom{1}}$  & $2.02\cdot 10^{-4}$  & $2.00\cdot 10^{-4}$ & $1.66$ & $1.67$ \\ 
    $2^{-4\phantom{1}}$  & $5.87\cdot 10^{-5}$  & $5.74\cdot 10^{-5}$ & $1.78$ & $1.81$ \\ 
    $2^{-5\phantom{1}}$  & $1.66\cdot 10^{-5}$  & $1.58\cdot 10^{-5}$ & $1.81$ & $1.86$ \\ 
    $2^{-6\phantom{1}}$  & $4.72\cdot 10^{-6}$  & $4.22\cdot 10^{-6}$ & $1.81$ & $1.91$ \\ 
    $2^{-7\phantom{1}}$  & $1.38\cdot 10^{-6}$  & $1.09\cdot 10^{-6}$ & $1.76$ & $1.95$ \\ 
    $2^{-8\phantom{1}}$  & $4.64\cdot 10^{-7}$  & $2.67\cdot 10^{-7}$ & $1.57$ & $2.03$ \\ 
    $2^{-9\phantom{1}}$ & $1.78\cdot 10^{-7}$  & $5.48\cdot 10^{-8}$ & $1.37$ & $2.28$ \\ 
    \hline  \hline
  \end{tabular}
  \caption{Experimental convergence of the MOL and the FTCS methods.}\label{tab:Exp_Convergence}
\end{center}
\end{table}

\section{Experimental Data Analysis}\label{sec: DATA_An}
As discussed earlier, the total amount of fluid absorbed by the specimen, defined in \eqref{eq: Q_defn}, represents an experimentally measurable quantity that establishes a connection between the mathematical model and physical experiments. Because of its intrinsic physical interpretation, $Q(t)$ is expected to exhibit a monotonically non-decreasing behavior until the saturation level is reached. Nevertheless, in certain instances, the experimental observations in \eqref{eq:Q_data} deviate from this expected monotonic trend, primarily due to measurement inaccuracies or instrumental limitations. This is particularly evident for natural materials. Moreover, for aerial mortars (GS, GSN and GSP), the filtration rate is exceedingly high, resulting in rapid saturation, which complicates the extraction of reliable information about the initial stages of the process. To address these challenges, here we propose a monotonicity-preserving data fitting and reconstruction methodology. This approach aims both to extract physically consistent insights from experiments and to generate synthetic data in scenarios where direct measurements are absent.

\subsection{Monotonicity-preserving differential reconstruction of data}\label{sub_sec: Data_Reconstruction}
Let $Q_{data}^i \approx Q(\tau_i),$ $i = 0, \dots, N_{data},$ represent the observed experimental values at time $\tau_i\in [0,\mathsf{T}].$ Denote by $\mathbb{D}$ and $\mathbb{L}^2$ the spaces of real differentiable and Lebesgue square integrable functions, respectively. Given the space of smooth monotone functions
\begin{equation}\label{eq: Functional_Space}
    \mathcal{F}=\left\{f(t): f\in\mathbb{D}, \ \log\left(\frac{d f}{d t}\right)\in\mathbb{D}, \ \dfrac{d}{dt} \! \left( \! \log\left(\frac{d f}{d t}\right) \! \right)\in \mathbb{L}^2 \right\},
\end{equation}
we aim to determine $Q^*\in\mathcal{F}$ which minimizes the mismatch with the measurements in \eqref{eq:Q_data}. Specifically, we seek for a solution to
\begin{equation}\label{eq:constrained_optim_problem}
    Q^*=\argmin_{f\in\mathcal{F}}\mathcal{L}(f,Q_{data}^0,\dots,Q_{data}^{N_{data}}),
\end{equation}
where $\mathcal{L}$  denotes a suitable loss function, whose explicit formulation will be detailed later on. The functional constrained optimization problem \eqref{eq:constrained_optim_problem} is particularly challenging to solve, even for relatively simple loss functions. 
Therefore, here we introduce a different approach built upon the following result, proved in \cite{Ramsay}, which provides a characterization for the space $\mathcal{F}.$ 
\begin{theorem}\label{thm: ODE_system}
    Consider the functional space $\mathcal{F}$ defined in \eqref{eq: Functional_Space}. Then, for each $f\in \mathcal{F},$ there exists a function $w_f\in \mathbb{L}^2$ such that $\bm{F}_w(t):=[f(t), \frac{df}{dt}(t)]^T$ solves the ODE system
    \begin{equation}\label{eq: ODE_System}
        \dfrac{d\bm{F}_w}{dt}(t)=\begin{bmatrix}
            0 & 1 \\ 0 & w_f(t)
        \end{bmatrix}\cdot \bm{F}_w(t).
    \end{equation}
\end{theorem}
The theoretical result established in Theorem \ref{thm: ODE_system} allows us to reformulate \eqref{eq:constrained_optim_problem} as the task of identifying the unconstrained function $w_f(t)$, represented as a linear combination of basis functions. Specifically, we employ the truncated Legendre expansion \cite{sym15061282} 
\begin{equation}\label{eq: Leg_expan}
    w^M_f(t) = \sum_{n=0}^{M} c_n \, P_n^{[0,\mathsf{T}]}(t) = \sum_{n=0}^{M} \frac{c_n}{n! \, \mathsf{T}^n} \frac{d^n}{dt^n} \big[(\mathsf{T} - t)^n t^n\big], \quad t \in [0, \mathsf{T}],
\end{equation}  
where $c_n$ are the unknown coefficients and $P_n^{[0,\mathsf{T}]}$, $n=0,\dots,M$, denote the shifted Legendre polynomials of degree $n,$ expressed through Rodrigues' formula (see, for instance, \cite{COHEN20121947,Askey2005}). Using \eqref{eq: ODE_System}, we then formulate the following Cauchy problem  
\begin{equation}\label{eq: Cauchy_Problem}
    \frac{d\bm{F}^M_w}{dt}(t) = \begin{bmatrix}
        0 & 1 \\ 
        0 & w^M_f(t)
    \end{bmatrix} \bm{F}^M_w(t), 
    \quad \bm{F}^M_w(0) = \begin{bmatrix}
        Q_{data}^0 \\ (Q_{data}^1 - Q_{data}^0)/(\tau_1 - \tau_0)
    \end{bmatrix},
\end{equation}  
whose solution $\bm{F}^M_w$ is expressed as a function of the expansion coefficients vector $\bm{c} = [c_0, \dots, c_M]^T \in \mathbb{R}^{M+1}$. Finally, we define the unconstrained minimization problem  
\begin{equation}\label{eq:UNconstrained_optim_problem}
    \bm{c}^* = \argmin_{\bm{c} \in \mathbb{R}^{M+1}} \left( \sum_{i=0}^{N_{data}} \left| Q_{data}^i - \bm{e}_1 \cdot \bm{F}^M_w(\tau_i; \bm{c}) \right|^2 + \lambda \|w^M_f(t; \bm{c})\|^2_{\mathbb{L}^2} \right),
\end{equation}  
where $\bm{e}_1 = [0, 1]$ is a row vector, $\lambda>0$ and 
\begin{equation*}
    \lambda\|w^M_f(t; \bm{c})\|^2_{\mathbb{L}^2} = \lambda\int_0^\mathsf{T} |w^M_f(t)|^2 \, dt = \lambda \mathsf{T} \sum_{n=0}^{M} \frac{c_n^2}{2n+1},
\end{equation*}  
is a regularization term obtained from \eqref{eq: Leg_expan} by exploiting the orthogonality of the shifted Legendre polynomials with respect to the Jacobi weight $W^{0,0}(x)=1.$

Here, we address problem \eqref{eq:UNconstrained_optim_problem} using the interior-point algorithm described in \cite{Byrd2000,Byrd99,Waltz2006} and implemented in Matlab's \texttt{fmincon} function. Specifically, we configure the solver with optimality and step tolerances of $10^{-14}$ and exploit parallelization for gradient computations. Furthermore, regarding the evaluation of the objective function in \eqref{eq:UNconstrained_optim_problem}, we approximate the solution to the ODE system \eqref{eq: Cauchy_Problem} using the Dormand-Prince explicit Runge-Kutta embedded method (cf. \cite{DORMAND198019} and \cite[Section 11.6]{Hairer_1987}) of Matlab's \texttt{ode45} routine (we refer to \cite{ODE_SUITE} for a comprehensive description of Matlab ODE Suite). 

\begin{figure}[htp]
     \centering
     \includegraphics[width=0.48\linewidth]{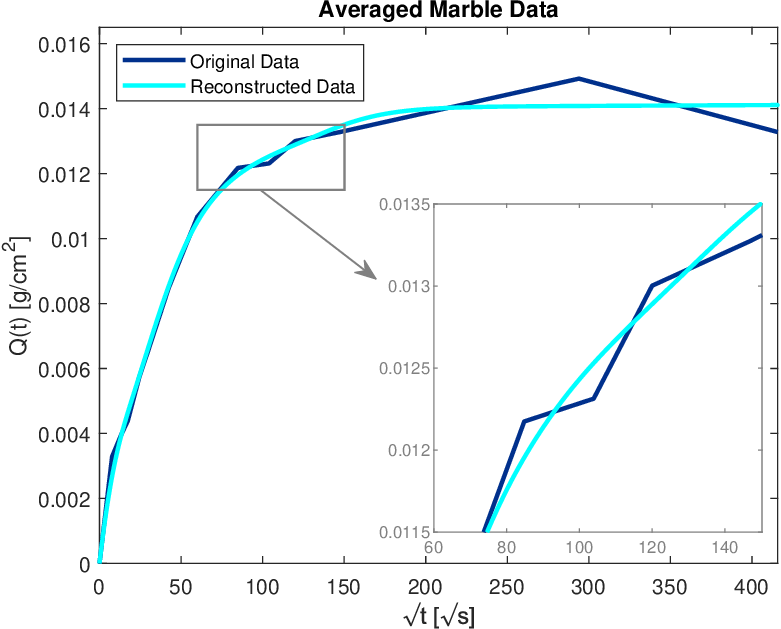} 
     \includegraphics[width=0.48\linewidth]{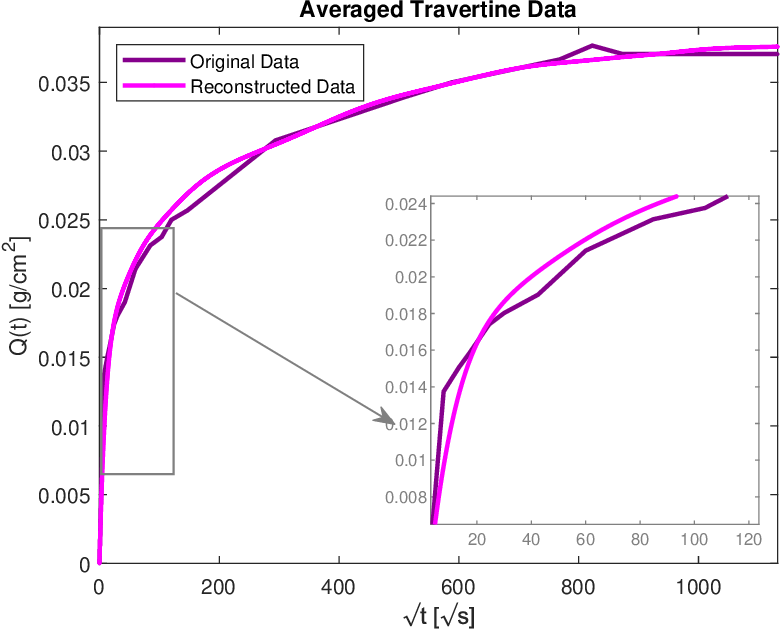}
     \caption{Monotonicity-preserving data differential reconstruction procedure applied to the averaged Carrara marble (left) and travertine (right) datasets.}
     \label{fig:Marble_Tr_Avg_Recon}
\end{figure}

\begin{figure}[htp]
    \centering
    \includegraphics[width=0.48\linewidth]{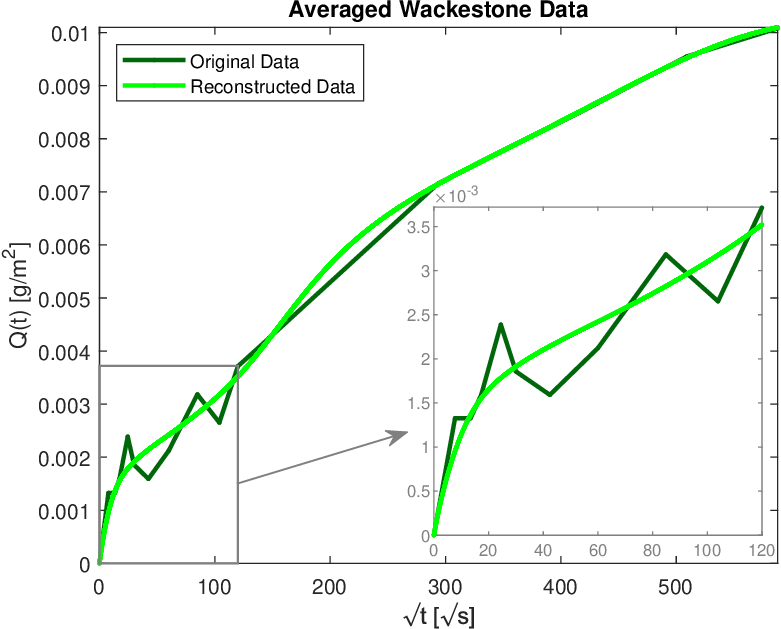}
    \includegraphics[width=0.48\linewidth]{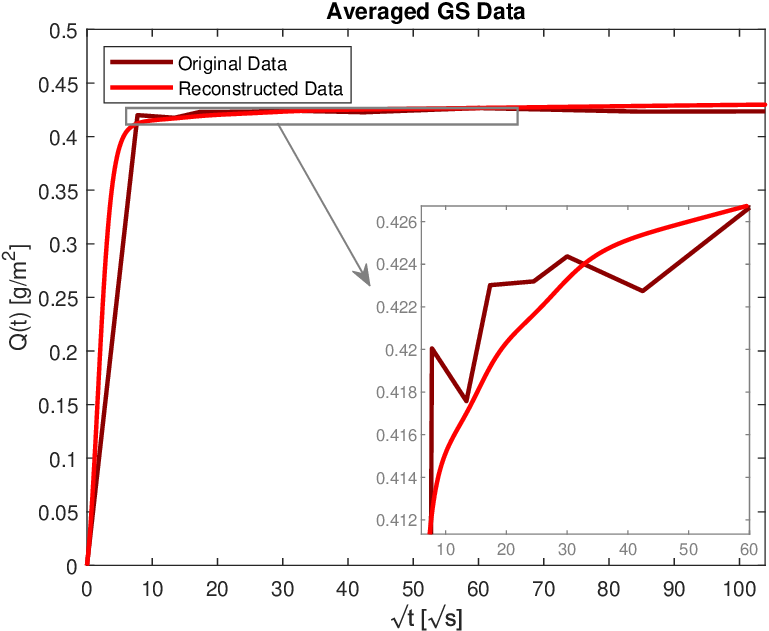}
    \caption{Monotonicity-preserving data differential reconstruction procedure applied to the averaged wackestone (left) and GS (right) datasets.}
    \label{fig:WS_GS_Recon}
\end{figure}

\begin{figure}[htp]
    \centering
    \includegraphics[width=0.48\linewidth]{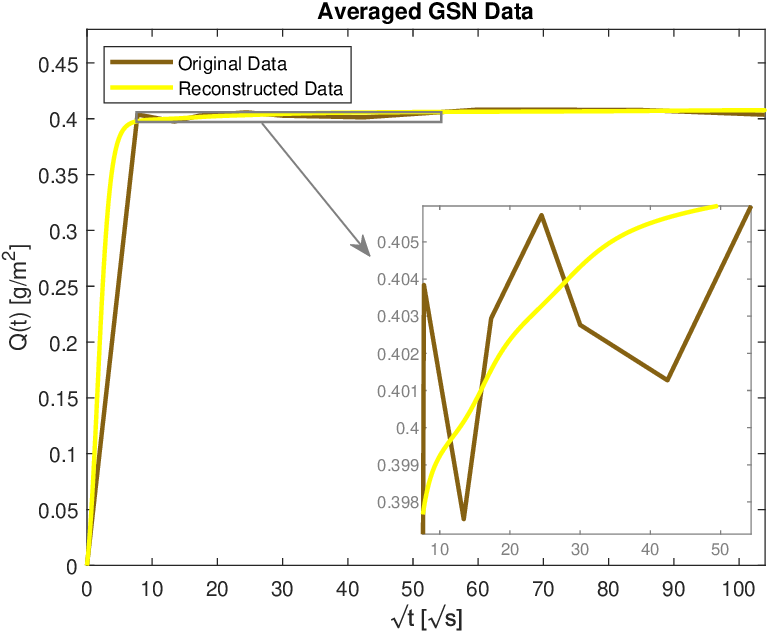}
    \includegraphics[width=0.48\linewidth]{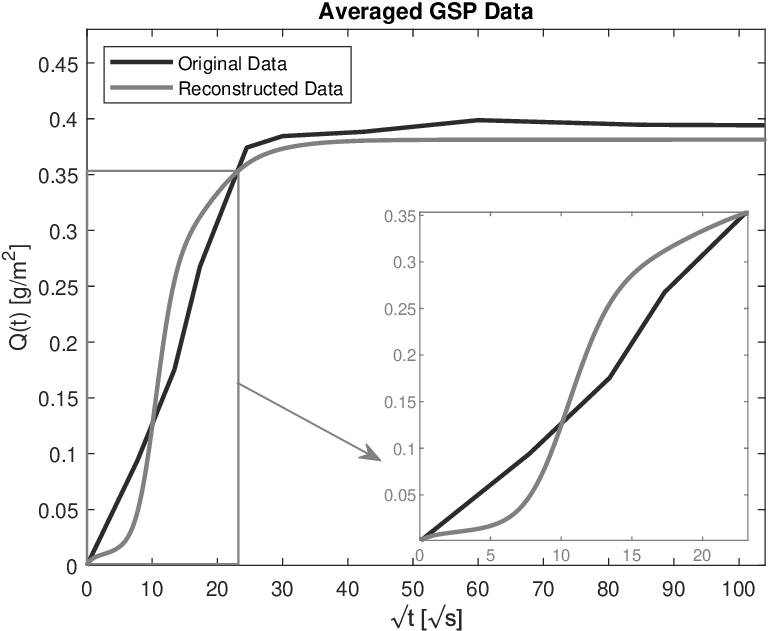}
    \caption{Monotonicity-preserving data differential reconstruction procedure applied to the averaged GSN (left) and GSP (right) datasets.}
    \label{fig:GSN_GSP_Recon}
\end{figure}

The proposed differential reconstruction procedure has been applied to the average datasets of all the materials, both natural and artificial, except for hydraulic mortars (OT), for which it is not necessary. The corresponding outcomes are reported in Figure \ref{fig:Marble_Tr_Avg_Recon} for Carrara marble and for the travertine, in Figure \ref{fig:WS_GS_Recon} for the wackestone and in Figures \ref{fig:WS_GS_Recon} and \ref{fig:GSN_GSP_Recon} for three kinds of aerial mortars. Furthermore, the dimensions of the polynomial space in which the expansion \eqref{eq: Leg_expan} is computed, as well as the values of the weight $\lambda$ for the different materials, are reported in Table \ref{table:lambda_values}. We remark that the value of $M$ has been chosen large enough to ensure a good representation capability and provide greater flexibility in the expansion, while still balancing the potential instability and numerical artifacts of the whole procedure.
\begin{table*}[!ht]
\centering
{
\begin{tabular}{|c|c|c|} \hline
\hline
 Material & $M$ & $\lambda$  \\\hline\hline
Carrara Marble & $25$& $3.000\cdot 10^{-6}$
\\
Tivoli Travertine & $25$& $5.200\cdot 10^{-4}$
\\
Wackestone & $25$ & $1.000\cdot 10^{-4}$\\
GS & $25$& $2.245\cdot 10^{-3}$ \\
GSN & $25$ & $2.000\cdot 10^{-3}$\\
GSP & $25$ & $2.000\cdot 10^{-3}$\\
\hline\hline
\end{tabular} }
\caption{Number of polynomials and regularization weights applied in the reconstruction of datasets for different materials.}\label{table:lambda_values}
\end{table*}

\begin{remark}
    The reconstruction procedure has been performed on the average datasets after applying appropriate nondimensionalization and linear rescaling to the data and variables onto $[0,1]$. Once completed, the reconstructed results are then mapped back to their original dimensional form by inverting these transformations.
\end{remark}

\section{Calibration of model parameters against imbibition data}\label{sec:Calibration}

As discussed in the previous sections, the determination of the function $B(\cdot)$ in \eqref{eq: B_def} that accurately replicates the capillary rise properties of the materials introduced in Section \ref{sec:material} represents a challenging task. Here, we adopt for $B(\cdot)$ the formulation \eqref{eq: B_shape} and consider the PDE model \eqref{eq:richards} complemented with the Robin boundary condition \eqref{eq:RBC}. Using the average experimental or reconstructed datasets, we then calibrate the following parameters:
\begin{itemize}
    \item $n_0\in\mathcal{N},$ the porosity of the material, accounting for the ranges in Table \ref{table:allmaterials},
    \item $s_{R}\in[0,1),$ the residual saturation of the porous medium,
    \item $s_{S}\in(0,1],$ the maximum saturation level of the porous medium,
    \item $D\in[0,\bar{D}]\subset\mathbb{R}^+,$ the diffusion coefficient of water within the porous medium, 
    \item $K_w\in[0,\bar{K}]\subset\mathbb{R}^+,$ the water exchange rate.
\end{itemize}

Define $\Omega= \mathcal{N} \times [0,1) \times (0,1] \times [0,\bar{D}] \times [0,\bar{K}] \subset \mathbb{R}^{5}$ and denote by $\bm{p} = [n_0,s_{R},s_{S},D,K_w]^T \in \Omega$ the vector of the parameters to be determined (with the values of $\bar{D}$ and $\bar{K}$ to be provided later). Let $\bm{Q}^{num}=\{Q^k(\bm{p}; \Delta z, \Delta t)\}_{0\leq k \leq N_t}$  denote the approximated total quantity of fluid absorbed by the specimen, computed using the MOL scheme \eqref{eq: Numerical_Scheme} as detailed in \eqref{eq: Q_num}. Assume that, for each $i=0,\dots,N_{data},$ there exists an integer $k_i\in\{0,\dots,N_t\}$ such that $\tau_i=k_i \Delta t,$ with $Q_{data}^i \approx Q(\tau_i).$ Then, the calibration procedure consists in finding a solution to
\begin{equation}\label{eq:Calibration_procedure}
    \bm{p}^* = \argmin_{\bm{p} \in \Omega} \left( \ell(\bm{Q}_{data},\bm{Q}^{num}(\bm{p};\Delta z,\Delta t)) + \varphi (Q_{data}^{N_{data}},Q^{k_{N_{data}}}(\bm{p};\Delta z,\Delta t) \right),
\end{equation}  
where the loss function $\ell(\cdot)$ measures the discrepancy between the observed and simulated dynamics throughout the entire time interval $[0,\mathsf{T}]$ and $\varphi (\cdot)$ represents the final cost that depends on the state of the system at the time $\mathsf{T}$. 

More specifically, denoted by $\mathcal{A}(\bm{p})=\mathcal{A}(\bm{Q}_{data},\bm{Q}^{num}(\bm{p};\Delta z,\Delta t))$ the set of all the possible warping paths $\pi=\{\{i,k_j\} \ | \ i,j\in [0,N_{data}]\cap \mathbb{N}\}$ between the experimental data and the simulation outcomes, we consider
\begin{equation}\label{eq:L2-DTW_Cost}
      \ell(\cdot)=\dfrac{\lambda_{2}}{N_{data}}\sum_{i=0}^{N_{data}} \dfrac{\left| Q_{data}^i - Q^{k_i}(\bm{p}; \Delta z, \Delta t) \right|^2}{(Q_{data}^i)^2} +\lambda_{DTW}\min_{\pi\in\mathcal{A}(\bm{p})} \Vert Q_{data}^i-Q^{k_j}(\bm{p}; \Delta z, \Delta t)\Vert_2,
\end{equation}
as a weighted sum of the Squared Relative Euclidean (SRE) and the Dynamic Time Warping (DTW)  errors (we refer to \cite{DTW1,DTW2} and references therein for further details) that quantify pointwise deviations and temporal misalignments of the data, respectively. Furthermore, we take the final cost as follows
\begin{equation}\label{eq:Final_Cost}
      \varphi (Q_{data}^{N_{data}},Q^{k_{N_{data}}}(\bm{p};\Delta z,\Delta t)=\left\lbrace\begin{array}{ll}\lambda_{\varphi}, \qquad & \text{if  } \ \vert Q_{data}^{N_{data}}-Q^{k_{N_{data}}}(\bm{p};\Delta z,\Delta t)\vert^{2}>\varepsilon_\varphi, \\ 
    0,& \text{otherwise,}
    \end{array}\right.
\end{equation}
with positive weight $\lambda_\varphi=10$ and threshold parameter $\varepsilon_\varphi=10^{-4}.$

\subsection{The Particle Swarm Multigrid Approach}\label{subsec:Particle_Swarm_MG}
The optimization problem \eqref{eq:Calibration_procedure}, arising from parameter estimation, is affected by the typical complications of inverse problems, such as strong nonlinearity, ill-posedness and high computational demands. Furthermore, the corresponding objective functional may be non-smooth and characterized by the presence of multiple local minima, which limits the effectiveness of gradient-based optimization methods. To address these difficulties, we adopt a gradient-free optimization strategy based on Particle Swarm Optimization (PSO) (see, for instance, \cite{PSO_overview,PSO}),  coupled with a multigrid scheme designed to enhance convergence and robustness across discretization levels, while simultaneously reducing computational costs. Regarding the PSO implementation, we employ the MATLAB built-in function \texttt{particleswarm} with the following settings: the swarm size is set to $1000$, the maximum number of iterations to $500$ and the maximum number of stall iterations to $50$. The function tolerance is set to $10^{-7}$, with both the self-adjustment and social-adjustment weights equal to $1.49$.

The proposed multigrid strategy relies on two resolution levels: a two-stages \hyperref[coarse_phase]{\emph{coarse phase}}, defined on a sparser spatio-temporal discretization, and a \hyperref[fine_phase]{\emph{fine phase}}, associated with a denser grid.

\begin{description}[topsep=\baselineskip, parsep=0pt]
    \item [Coarse phase]\label{coarse_phase} 
    In the first stage, two independent optimization subproblems are solved in parallel with PSO. Each subproblem involves a simplified version of the objective functional \eqref{eq:Calibration_procedure}, where only one of the two contributions appearing in \eqref{eq:L2-DTW_Cost} is retained. The resulting errors, here referred to as $\varepsilon_2$ and $\varepsilon_{DTW},$ are then employed to compute the weights $\lambda_{2} =1/\varepsilon_2$ and $\lambda_{DTW} = 1/\varepsilon_{DTW}.$ These values, computed for each material, are reported in Table \ref{tab:weights_Cal_values}. \\
    In the second stage, the PSO is employed to minimize, on the coarse grid, the full functional \eqref{eq:Calibration_procedure}, equipped with the weights retrieved from the previous step. For the entire coarse phase, the admissible porosity ranges $\mathcal{N}$ are set as reported in Table \ref{table:allmaterials}. In compliance with the experimental findings of \cite{saturation_range}, we consider $s_R\in[0.1, 0.75]$ and $s_S\in[0.5, 0.98].$ Furthermore, we set $\bar{D}=0.1$ and take $\bar{K}=100$ large enough to investigate both Neumann-like (lower values of $K_w$) and Dirichlet-like (higher values of $K_w$) scenarios. Therefore, the parameter space for both the coarse stages reads $\Omega= \mathcal{N} \times [0.1,0.75] \times [0.5,0.98] \times [0,0.1] \times [0,100] \subset \mathbb{R}^{5}.$ \\
    
    \item [Fine phase]\label{fine_phase}
    This phase consists of a sequence of $\nu \in \mathbb{N}$ optimization steps, each involving the minimization of the full objective functional \eqref{eq:Calibration_procedure} over a progressively refined parameter space $\Omega^{(n)}$, with $n = 1,\dots,\nu$. Let $\bm{p}^{*,0}\in \Omega$ represent the optimal solution retrieved from the second step of the \hyperref[coarse_phase]{\emph{coarse} phase}. Each subsequent space $\Omega^{(n)}$ is recursively defined around the optimum $\bm{p}^{*(n-1)}$ on $\Omega^{(n-1)}$ as
    \begin{equation*}\label{eq:refined_param_space}
        \Omega^{(n)} = \prod_{i=1}^5 \left[p_i^{*(n-1)} (1- \sigma_i^{(n)}),\, p_i^{*(n-1)} (1+ \sigma_i^{(n)})\right], \qquad \qquad n = 1,\dots,\nu,
    \end{equation*}
    where $p_i^{*(n-1)},$ $i=1,\dots,5,$ denotes the $i$-th component of $\bm{p}^{*(n-1)}$ and $\sigma_i^{(n)}\in(0,1)$ is a user-defined local refinement radius (here, we set $\sigma_i^{(n)}=(2n)^{-1}$ for each $i$). Furthermore, for each optimization step of this phase, the values of the stepsizes in \eqref{eq: Numerical_Scheme} are chosen to satisfy the CFL condition \eqref{eq: Stability_Condition} at the least refined level.
    \end{description}
    
\begin{table}[htbp]
    \centering
    {\begin{tabular}{|c|c|c|}
        \hline\hline
        Material & $\lambda_{2}$ & $\lambda_{DTW}$ \\
        \hline \hline
        Carrara Marble & $1.20 \cdot 10^{-3}$ & $1.65\cdot 10^{-3}$ \\
        Tivoli Travertine & $3.34\cdot 10^{-2}$ & $9.13\cdot 10^{-2}$ \\
        Travertine (parallel) & $6.41\cdot 10^{-2}$ & $1.19\cdot 10^{-1}$\\
        Travertine (perpendicular) & $7.35\cdot 10^{-3}$ &$4.41\cdot 10^{-1}$\\
        Wackestone & $5.75\cdot 10^{-1}$ & $1.35 \cdot 10^{0\phantom{-}}$ \\
        OT1  & $1.87 \cdot 10^{0\phantom{-}}$ & $1.77\cdot 10^{-1}$ \\
        OT2  & $6.61\cdot 10^{-1}$ & $9.50\cdot 10^{-3}$ \\
        OT3  & $1.40\cdot 10^{-1}$ & $4.70\cdot 10^{-3}$ \\
        OT4  & $1.90\cdot 10^{-2}$ & $4.61\cdot 10^{-3}$ \\
        GS  & $1.22 \cdot 10^{0\phantom{-}}$ & $6.91\cdot 10^{-3}$ \\
        GSN  & $9.52 \cdot 10^{0\phantom{-}}$ & $5.89\cdot 10^{-3}$ \\
        GSP  & $4.44\cdot 10^{-1}$ & $1.73 \cdot 10^{0\phantom{-}}$ \\
        \hline\hline
    \end{tabular}}
    \caption{Optimal weights for the loss function contributions associated with different materials, as determined during the \hyperref[coarse_phase]{coarse phase}.}
    \label{tab:weights_Cal_values}
\end{table}

\subsection{Calibration results for different materials}\label{subsec:Calibration_Results}

\begin{figure}[!ht]
    \centering
    \includegraphics[width=0.48\linewidth]{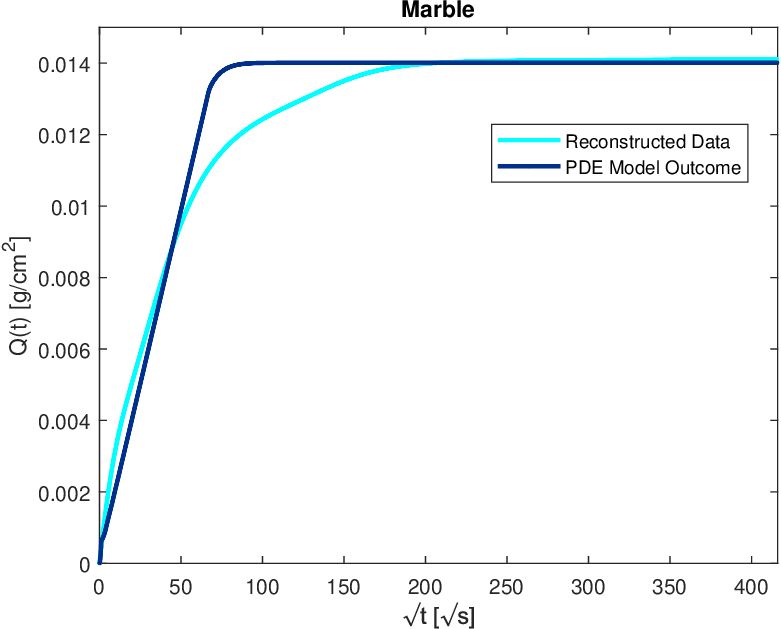}
    \includegraphics[width=0.48\linewidth]{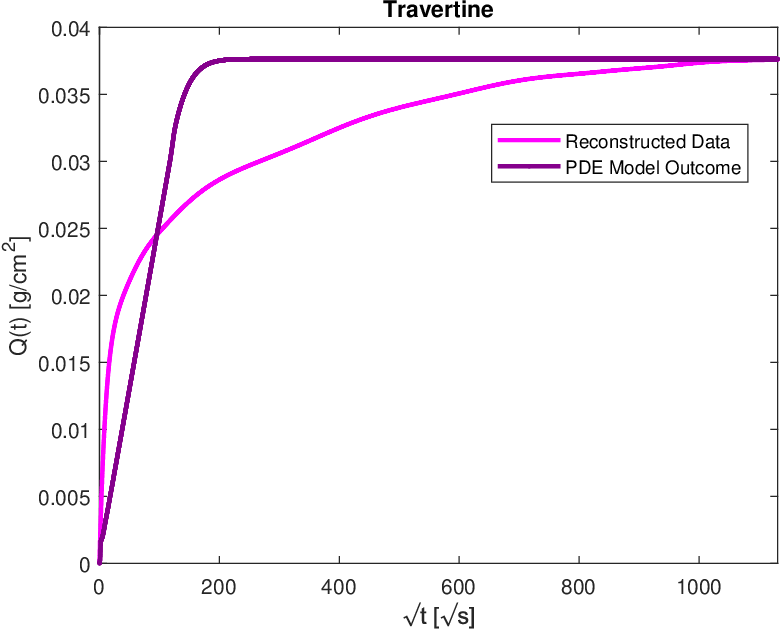}
    \caption{Simulation of the data-informed mathematical model, compared against the reconstructed experimental imbibition curve for Carrara marble (left) and travertine (right).}
    \label{fig:Marble_Travertine_Outcomes}
\end{figure}

\begin{figure}[!ht]
    \centering
    \includegraphics[width=0.48\linewidth]{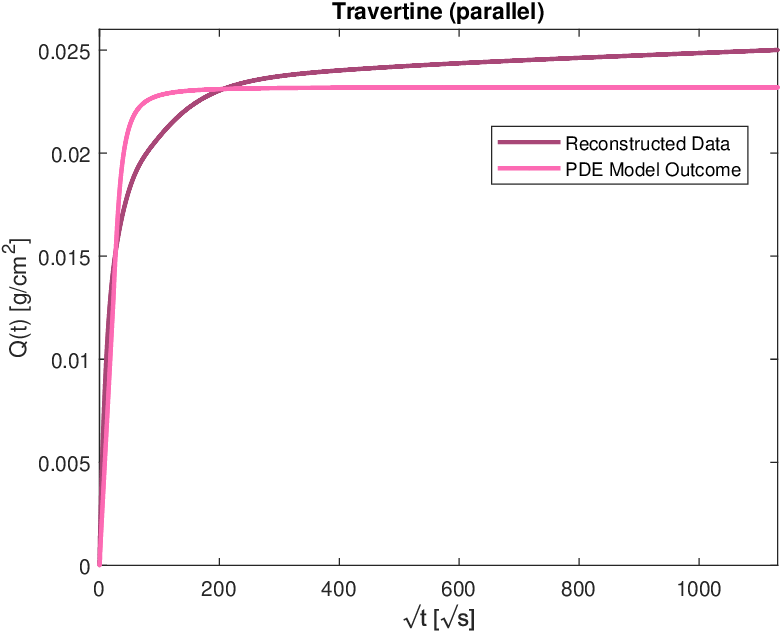}
    \includegraphics[width=0.48\linewidth]{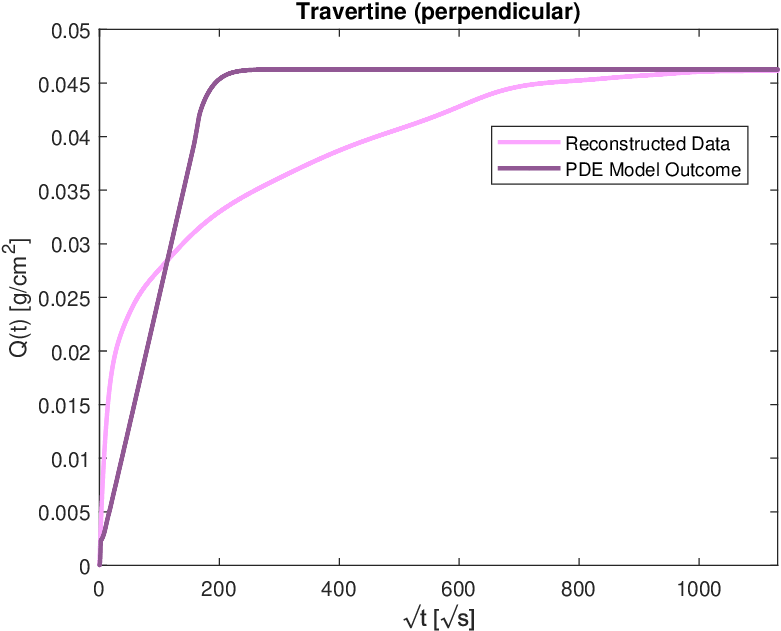}
    \caption{Simulation of the data-informed mathematical model, compared against the reconstructed experimental imbibition curves for Tivoli travertine with parallel (left) and perpendicular (right) anisotropy planes.}
    \label{fig:Anisotropic_Outcomes}
\end{figure}

\begin{figure}[!ht]
    \centering
    \includegraphics[width=0.48\linewidth]{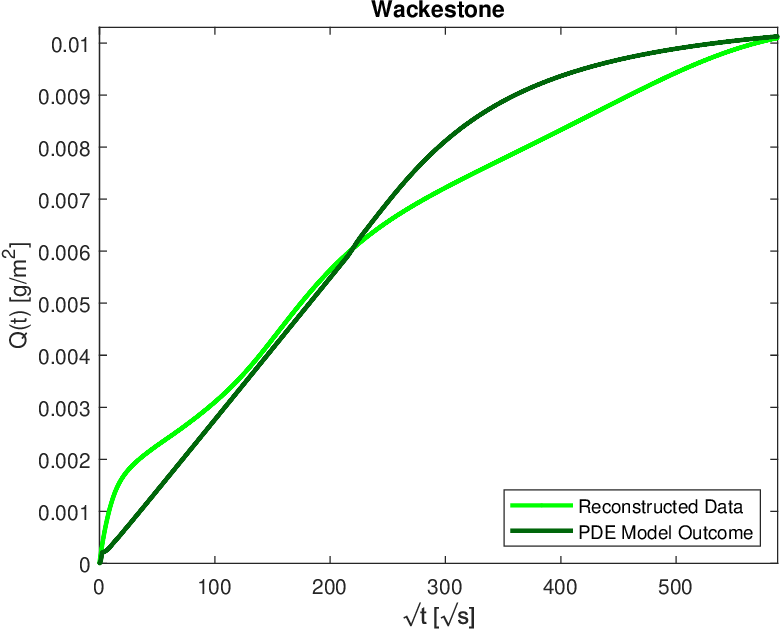}
    \includegraphics[width=0.48\linewidth]{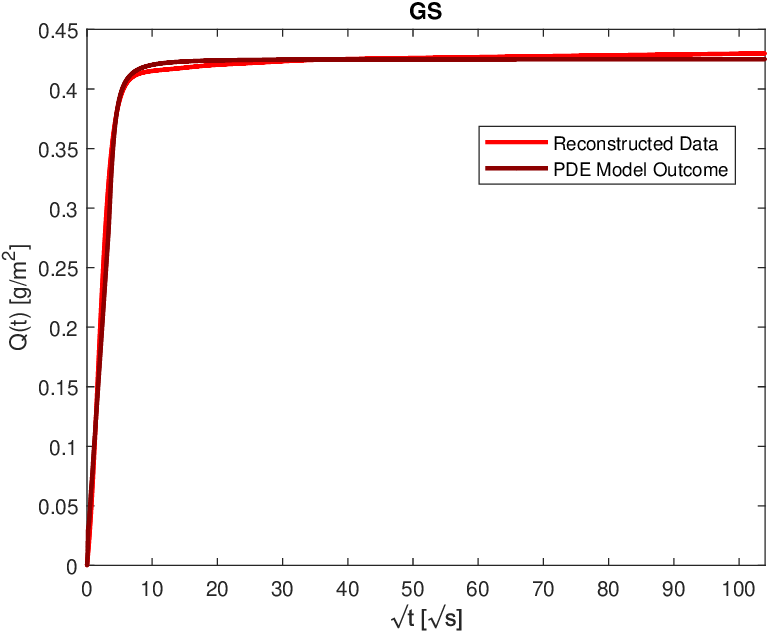}
    \caption{Simulation of the data-informed mathematical model, compared against the reconstructed experimental imbibition curves for wackestone (left) and GS (right).}
    \label{fig:Wackestone_GS_Outcomes}
\end{figure}

\begin{figure}[htp]
    \centering
    \includegraphics[width=0.48\linewidth]{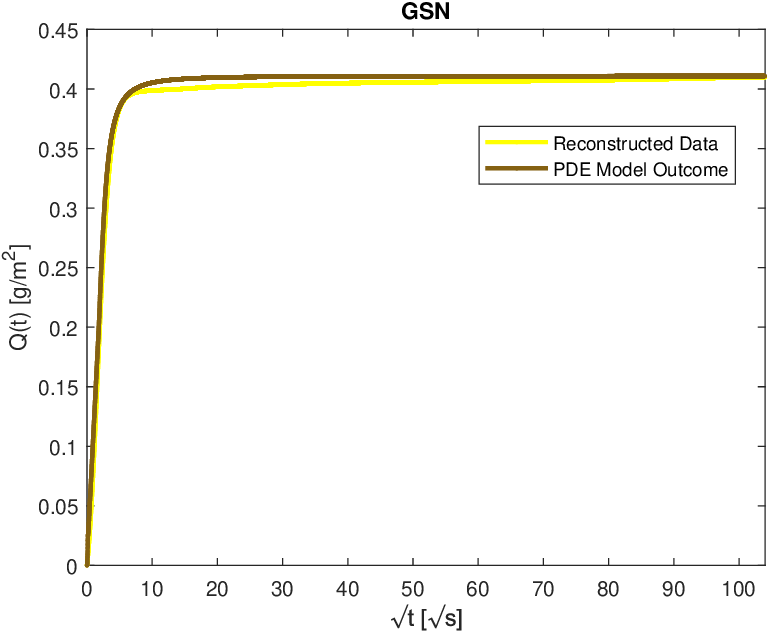}
    \includegraphics[width=0.48\linewidth]{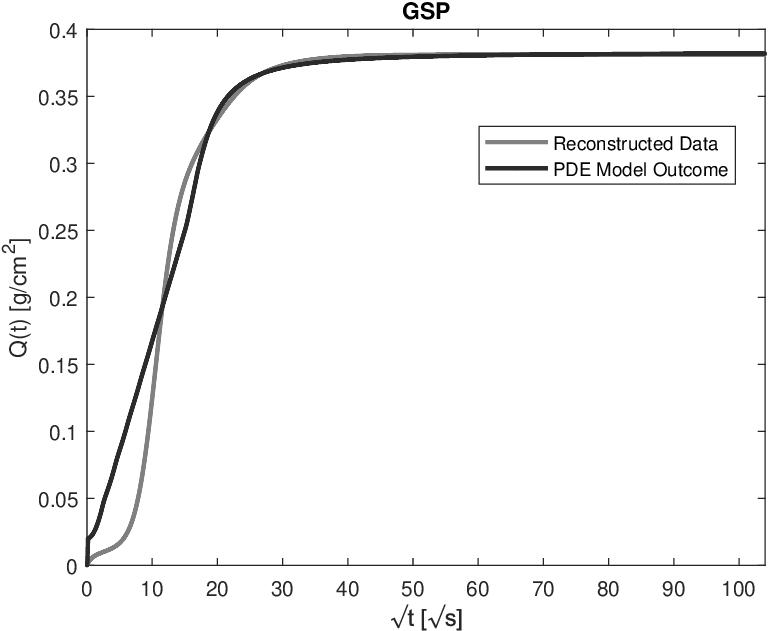}
    \caption{Simulation of the data-informed mathematical model, compared against the reconstructed experimental imbibition curves for GSN (left) and GSP (right).}
    \label{fig:GSN_GSP_Outcomes}
\end{figure}

\begin{figure}
\centering
    \includegraphics[width=0.48\linewidth]{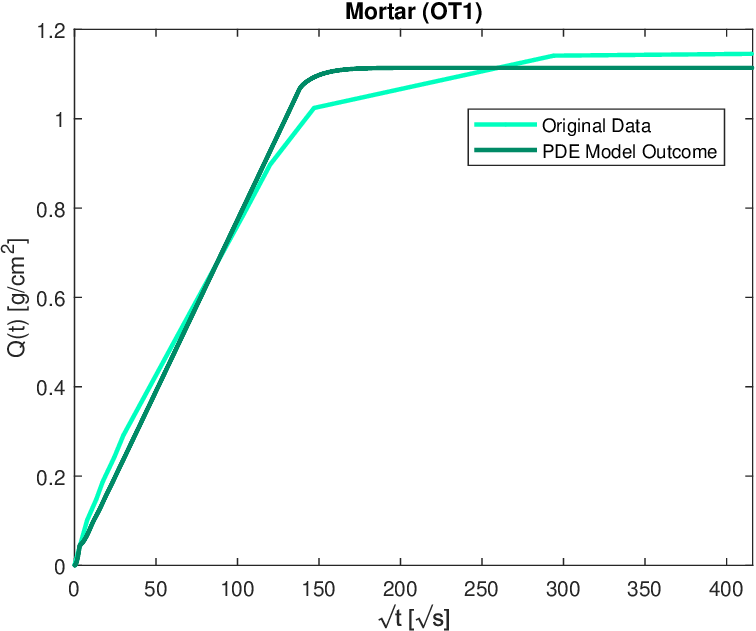} 
    \includegraphics[width=0.48\linewidth]{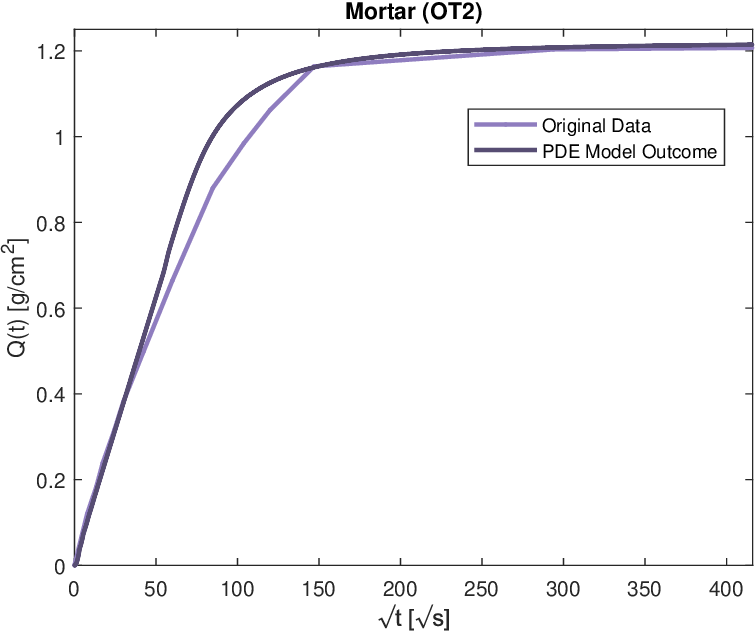}  \\
    \includegraphics[width=0.48\linewidth]{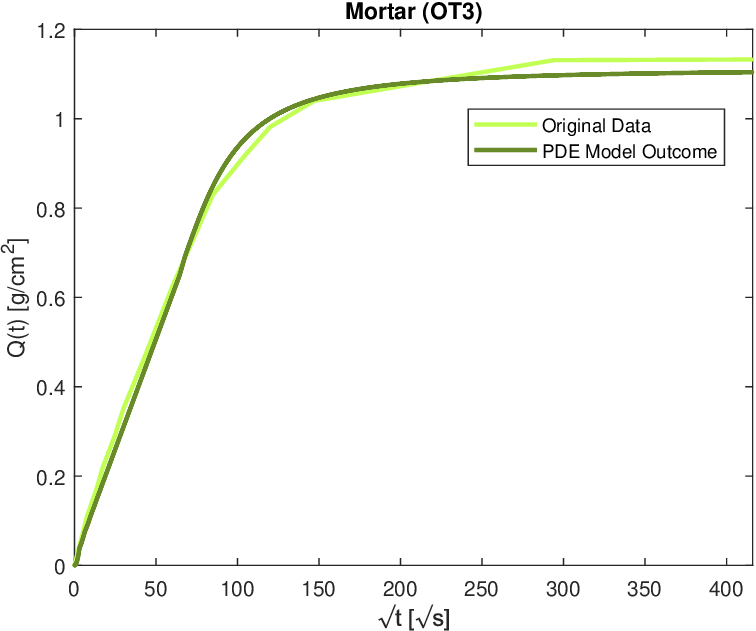}
    \includegraphics[width=0.48\linewidth]{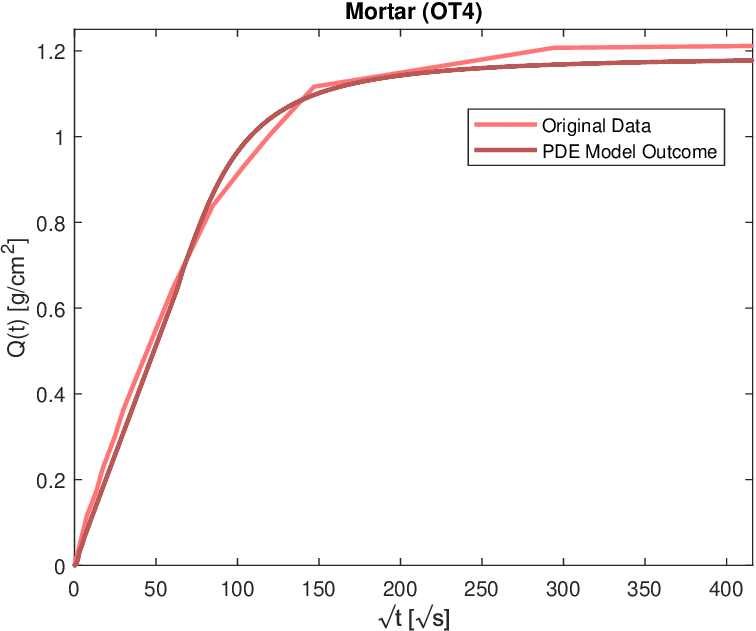}
    \caption{Simulation of the data-informed mathematical model, compared against the  experimental imbibition curves for different kinds of mortars.}
    \label{fig:OT_Outcomes}
\end{figure}

The PSO–multigrid calibration procedure has been applied to all the materials described in Section \ref{sec:material}. The corresponding optimized parameters and error metrics are reported in Tables \ref{tab:Natural_Materials_Calibration} and \ref{tab:Anisotropic} for natural materials and in Tables \ref{tab:GS_Calibr} and \ref{tab:OT_Calibr} for artificial ones. A direct comparison between the reconstructed data and the simulation results with calibrated parameters is provided with the Figures \ref{fig:Marble_Travertine_Outcomes}--\ref{fig:OT_Outcomes}. From a qualitative standpoint, the optimal parameter configurations generally yield simulated imbibition curves that adhere to the experimental profiles, effectively reproducing the dynamics of the observed phenomena.

A notable exception is represented by travertine. For this material, the discrepancy can plausibly be attributed to the presence of anisotropy planes in the material structure, features that may be present in the tested specimens but are not accounted for in the PDE model \eqref{eq:richards}. In addition, the experimental data reported in Figures \ref{fig:Marble_Tr_Avg_Recon} and \ref{fig:Marble_Travertine_Outcomes} (right panels) may reflect an average behavior across specimens characterized by significant internal heterogeneity, a typical consequence of the complex geological processes involved in travertine formation. This heterogeneity may manifest in the orientation of anisotropy planes, which can vary from one specimen to another. For this reason, we introduce a distinction between two classes of samples: those with anisotropy planes oriented parallel to the water flow in the specimen, and those in which the planes are oriented perpendicularly. The results of the parameter estimation for the two classes of travertine are reported in Table \ref{tab:Anisotropic} and further illustrated in Figure \ref{fig:Anisotropic_Outcomes}. This classification provides a more refined interpretation of the experimental outcomes and highlights the relevance of structural anisotropy in the observed deviations.

In general, from a quantitative point of view, the final calibration error, defined as the value of the objective functional in \eqref{eq:Calibration_procedure}, remains low for all the materials, including the travertine one.  This outcome highlights the effectiveness of the combined use of the SRE and the DTW errors within the cost functional, which allows for a balanced evaluation of both amplitude and shape differences in the calibration process.
\begin{table}[htbp]
\begin{center}
  \begin{tabular}{|c|c|c|c|} 
  \hline
  \multicolumn{4}{|c|}{Natural Materials}  \\
  \hline  \hline
   Parameters &  Marble & Travertine& Wackestone\\
   \hline \hline
    Residual saturation $s_{R}$ & $2.27\cdot 10^{-1}$ & $2.00\cdot 10^{-1}$&$2.40\cdot 10^{-1}$\\
    Maximum saturation $s_{S}$ &  $8.84\cdot 10^{-1}$ & $8.16\cdot 10^{-1}$&$9.50\cdot 10^{-1}$\\
    Diffusion rate $D$ & $1.09\cdot 10^{-5} $&$9.99\cdot 10^{-6}$&$5.70\cdot 10^{-7}$\\
        Water exchange rate $K_{w}$ &$8.57 \cdot 10^1$&$1.00 \cdot 10^2$&$0$ \\
    Porosity $n_{0}$ & $0.63\%$ &$1.33\%$&$0.22\%$\\
    SRE error & $2.94\cdot 10^{-2}$ & $1.25\cdot 10^{-1}$ & $1.29\cdot 10^{-1}$ \\
    DTW error &  $3.79\cdot 10^{-2}$ &$2.85\cdot 10^{-2}$& $3.70\cdot 10^{-2}$\\
    Objective function value & $1.04\cdot 10^{-4}$&$7.40\cdot 10^{-3}$ & $2.34\cdot 10^{-2}$ \\
    \hline\hline
  \end{tabular}
  \caption{Optimal parameters for the natural materials. 
  }\label{tab:Natural_Materials_Calibration}
\end{center}
\end{table}
\begin{table}[htbp]
\begin{center}
  \begin{tabular}{|c|c|c|} 
  \hline
  \multicolumn{3}{|c|}{Anisotropic Natural Materials}  \\
  \hline  \hline
   Parameters &  Travertine (parallel)& Travertine (perpendicular)\\
   \hline \hline
    Residual saturation $s_{R}$ & $2.40\cdot 10^{-1}$ & $2.22\cdot 10^{-1}$\\
    Maximum saturation $s_{S}$ &  $9.26\cdot 10^{-1}$ & $6.94\cdot 10^{-1}$\\
    Diffusion rate $D$ & $1.00\cdot 10^{-4} $&$9.19\cdot 10^{-6}$\\
        Water exchange rate $K_{w}$ &$8.14 \cdot 10^{1\phantom{-}}$&$7.14 \cdot 10^{1\phantom{-}}$ \\
    Porosity $n_{0}$ & $0.50\%$ &$1.86\%$\\
    SRE error & $2.60\cdot 10^{-2}$ & $1.29\cdot 10^{-1}$ \\
    DTW error &  $8.50\cdot 10^{-3}$ &$3.29\cdot 10^{-2}$\\
    Objective function value & $6.04\cdot 10^{-4}$&$8.40\cdot 10^{-3}$ \\
    \hline\hline
  \end{tabular}
  \caption{Optimal parameters for travertine groups classified by anisotropy plane orientation.
  }\label{tab:Anisotropic}
\end{center}
\end{table}
\begin{table}[htbp]
\centering
\begin{tabular}{|c|c|c|c|}
\hline
  \multicolumn{4}{|c|}{Artificial Materials}  \\
\hline  \hline
Parameters & GSN & GS & GSP \\\hline\hline
Residual saturation $s_{R}$ & $2.00\cdot 10^{-1}$ & $3.77\cdot 10^{-1}$ & $5.34\cdot 10^{-1}$ \\
Maximum saturation $s_{S}$ & $9.44\cdot 10^{-1}$ & $7.02\cdot 10^{-1}$ & $9.50\cdot 10^{-1}$ \\
Diffusion rate $D$ & $7.55\cdot 10^{-2}$ & $8.60\cdot 10^{-2}$ & $2.75\cdot 10^{-2}$ \\
Water exchange rate $K_{w}$ & $0$ & $0$ & $0$ \\
Porosity $n_{0}$ & $24.20\%$ & $29.61\%$ & $20.09\%$ \\
SRE error & $2.10\cdot 10^{-4}$ & $7.15\cdot 10^{-5}$ & $3.30\cdot 10^{-1}$ \\
DTW error & $1.75\cdot 10^{-2}$ & $1.10\cdot 10^{-2}$ & $8.00\cdot 10^{-2}$ \\
Objective function value & $6.00\cdot 10^{-4}$ & $3.75\cdot 10^{-6}$ & $8.10\cdot 10^{-3}$ \\
\hline\hline
\end{tabular}
\caption{Optimal parameters for the artificial materials.}
\label{tab:GS_Calibr}
\end{table}
\begin{table}[htbp]
\centering
\begin{tabular}{|c|c|c|c|c|}
\hline
  \multicolumn{5}{|c|}{Artificial Materials}  \\
\hline  \hline
Parameters & OT1 & OT2 & OT3 & OT4 \\\hline\hline
Residual saturation $s_{R}$ & $3.19\cdot 10^{-1}$ & $2.00\cdot 10^{-1}$ & $2.00\cdot 10^{-1}$ & $2.00\cdot 10^{-1}$ \\
Maximum saturation $s_{S}$ & $9.50\cdot 10^{-1}$ & $7.81\cdot 10^{-1}$ & $6.83\cdot 10^{-1}$ & $8.80\cdot 10^{-1}$ \\
Diffusion rate $D$ & $2.23\cdot 10^{-4}$ & $9.00\cdot 10^{-4}$ & $7.51\cdot 10^{-4}$ & $5.40\cdot 10^{-4}$ \\
Water exchange rate $K_{w}$ & $1.00 \cdot 10^2$ & $0$ & $0$ & $6.00\cdot 10^{-4}$ \\
Porosity $n_{0}$ & $44.33\%$ & $38.80\%$ & $40.20\%$ & $33.50\%$ \\
SRE error & $1.80\cdot 10^{-2}$ & $4.00\cdot 10^{-3}$ & $6.70\cdot 10^{-3}$ & $1.87\cdot 10^{-2}$ \\
DTW error & $1.34\cdot 10^{-1}$ & $2.23\cdot 10^{-1}$ & $9.75\cdot 10^{-2}$ & $1.39\cdot 10^{-1}$ \\
Objective function value & $9.85\cdot 10^{-3}$ & $5.32\cdot 10^{-4}$ & $4.75\cdot 10^{-4}$ & $1.00\cdot 10^{-3}$ \\
\hline\hline
\end{tabular}
\caption{Optimal parameters for the OT samples. Here, the values of $n_0$ have been provided by the producer.}
\label{tab:OT_Calibr}
\end{table}
\section{Conclusions and future perspectives}\label{sec:concl}

This work presents a comprehensive, data-informed mathematical framework specifically designed to characterize capillary water absorption in porous materials commonly employed in contexts of relevance to heritage science. The proposed approach integrates targeted laboratory experimentation with well-established mathematical models to yield a material-specific quantification of key physical parameters. Specifically, the imbibition dynamics is described through Richards-type differential equation built upon Darcy’s law, and simulated via a second order explicit numerical scheme based on the method of lines technique. On the experimental side, a series of imbibition experiments have been performed on various natural and artificial porous materials, generating a broad dataset. Additionally, a novel monotonicity-preserving reconstruction technique has been developed to correct slight irregularities in the experiments data, attributed to measurement noise and instrumental limitations, while retaining the physical coherence of the underlying monotone water uptake process. A dedicated calibration procedure, based on particle swarm optimization and enhanced by a multigrid hierarchical strategy, is introduced for the data-driven identification of the model parameters. The integration of the squared relative error, of the dynamic time warping distance and of an endpoint-discrepancy penalization term within the objective functional has effectively ensured accurate capture of both the overall evolution and the final state of the imbibition process. In conclusion, the outcomes for the monitored materials confirm the reliability of the comprehensive experimental-mathematical methodology we designed, which proves to be a versatile tool for analyzing water absorption phenomena in porous geomaterials, with potential applications in cultural heritage preservation.

Given the promising results of this study, several potential directions for future work emerge. A natural extension would involve incorporating anisotropy and heterogeneous permeability to better represent complex materials like travertine. Further investigation into additional degradation mechanisms, such as salt crystallization and pollutant transport, is also of high interest. Finally, applying the proposed techniques to in-situ datasets from historical structures could enable the model’s deployment as a diagnostic and predictive tool for real-world conservation efforts.
 

\section*{Acknowledgements.} G. B. and M. P. are members of the Gruppo Nazionale Calcolo Scientifico-Istituto Nazionale di Alta Matematica (GNCS-INdAM). G. B. and E. C. B. are in the PRIN project MATHPROCULT Prot. P20228HZWR, CUP B53D23015940001.  L. M., M. D. F. and M. P. are involved in the Project PE0000020 CHANGES - CUP IAC-CNR: B53C22003890006, CUP Sapienza:  B53C22003780006, PNRR Mission 4 Component 2 Investment 1.3, Funded by the European Union - NextGenerationEU" under the Italian Ministry of University and Research (MUR).

\section*{Data availability statement}
The data that support the findings of this study are available from the corresponding author, M.P., upon reasonable request. 

\bibliographystyle{plainurl}
\bibliography{reference}

\end{document}